\documentclass[10pt,a4paper]{article}
\usepackage[left=1in, right=1in, top=1in, bottom=1in]{geometry}

\usepackage[utf8]{inputenc}
\usepackage{amsmath}
\usepackage{mathtools}
\usepackage{amssymb}
\usepackage{amsthm}
\usepackage{color}
\usepackage{enumitem}
\usepackage{tikz}
\usepackage{pgfplots}
\pgfplotsset{compat=newest}
\usetikzlibrary{plotmarks}
\usetikzlibrary{arrows.meta}
\usepackage{xcolor}
\usepackage{esint}
\usepackage{array}
\usepackage{url}
\usepackage{subfig}
\usepackage{graphicx}

\DeclareMathOperator{\R}{\mathbb{R}}

\DeclareMathOperator{\Pp}{\mathbb{P}}

\newcommand{\Ccinf}[1]{C^{\infty}_c\left(#1\right)}

\newcommand{\E}[1]{\ensuremath{\mathbb{E}\left[#1\right]}}
\newcommand{\var}[1]{\ensuremath{\mathrm{Var}\left(#1\right)}}

\newcommand{\norm}[1]{\ensuremath{\left\lVert#1\right\rVert}}
\newcommand{\abs}[1]{\ensuremath{\left\lvert#1\right\rvert}}
\newcommand{\set}[1]{\ensuremath{\left\lbrace#1\right\rbrace}}

\RequirePackage[capitalize,nameinlink]{cleveref}[0.19]
\crefname{section}{section}{sections}
\crefname{subsection}{subsection}{subsections}
\Crefname{corollary}{Corollary}{Corollaries}
\Crefname{proposition}{Proposition}{Propositions}
\Crefname{lemma}{Lemma}{Lemmas}
\Crefname{figure}{Figure}{Figures}
\crefformat{equation}{\textup{#2(#1)#3}}
\crefrangeformat{equation}{\textup{#3(#1)#4--#5(#2)#6}}
\crefmultiformat{equation}{\textup{#2(#1)#3}}{ and \textup{#2(#1)#3}}
{, \textup{#2(#1)#3}}{, and \textup{#2(#1)#3}}
\crefrangemultiformat{equation}{\textup{#3(#1)#4--#5(#2)#6}}%
{ and \textup{#3(#1)#4--#5(#2)#6}}{, \textup{#3(#1)#4--#5(#2)#6}}{, and \textup{#3(#1)#4--#5(#2)#6}}
\Crefformat{equation}{#2Equation~\textup{(#1)}#3}
\Crefrangeformat{equation}{Equations~\textup{#3(#1)#4--#5(#2)#6}}
\Crefmultiformat{equation}{Equations~\textup{#2(#1)#3}}{ and \textup{#2(#1)#3}}
{, \textup{#2(#1)#3}}{, and \textup{#2(#1)#3}}
\Crefrangemultiformat{equation}{Equations~\textup{#3(#1)#4--#5(#2)#6}}%
{ and \textup{#3(#1)#4--#5(#2)#6}}{, \textup{#3(#1)#4--#5(#2)#6}}{, and \textup{#3(#1)#4--#5(#2)#6}}
\crefdefaultlabelformat{#2\textup{#1}#3}

\newtheorem{theorem}{Theorem}
\newtheorem{corollary}[theorem]{Corollary}
\newtheorem{lemma}[theorem]{Lemma}
\newtheorem{proposition}[theorem]{Proposition}
\newtheorem{definition}[theorem]{Definition}

\newtheorem{remark}{Remark}

\numberwithin{theorem}{section}

\allowdisplaybreaks

\newcommand{\imgdir}{./img}

\title{Analytical and numerical study of a modified cell problem for the numerical homogenization of multiscale random fields.}
\author{Abdulle, A., Arjmand, D. and Paganoni, E.}

\begin{document}
\maketitle

\begin{abstract}
A central question in numerical homogenization of partial differential equations with multiscale coefficients is the accurate computation of effective quantities, such as the homogenized coefficients. Computing homogenized coefficients requires solving local corrector problems followed by upscaling relevant local data. The most naive way of computing homogenized coefficients is by solving a local elliptic problem, which is known to suffer from the so-called resonance error dominating all other errors inherent in multiscale computations. A far more efficient modelling strategy, based on adding an exponential correction term to the standard local elliptic problem, has recently been proved to result in exponentially decaying error bounds with respect to the size of the local geometry. The questions in relation with the accuracy and computational efficiency of this approach has been previously addressed in the context of periodic homogenization. The present article concerns the extension of mathematical and numerical study of this modified elliptic corrector problem to stochastic homogenization problems. In particular, we assume a stationary, ergodic micro-structure and i) establish the well-posedness of the corrector equation, ii) analyse the bias (or the systematic error) originating from additional exponential correction term in the model. Numerical results corroborating our theoretical findings are presented. 
\end{abstract}



\section{Introduction}
\label{sec: intro}	
A common denominator of many problems in natural sciences and engineering is that they include elements of heterogeneity and randomness at fine scales. For instance, the characteristics of widely used construction materials such as clay, ceramics, and concrete, or properties of polymer reinforced composite materials used in aerospace industry, are not only microscopically inhomogeneous, but also incorporate uncertainty due to imperfections distributed randomly in the microstructure \cite{Hil63,KFG03}. Other related applications include but are not limited to problems in solid and fluid mechanics, hydrogeology and chemistry, and transport of pollutant in porous media \cite{BeB91,BDI19}. Partial differential equations (PDEs) with multiscale random input data can be employed to model such  physical phenomena. In general, microscopic features influence the overall macroscopic behaviour, and therefore their effects have to be accurately quantified. From the point of view the numerical simulations, the problem is twofold: the presence of rapidly varying model parameters require resolutions
down to the fine scales, and many replicas of an already expensive multiscale problem/solution are needed to capture statistically meaningful quantities.

In order  to bypass the above mentioned computational challenges, one may resort to analytical tools, such as the homogenization theory, see e.g., \cite{BLP78,CDG02,JKO94},  which aims to describe the macroscopic behaviour of a microscopically inhomogeneous system. Once the effective quantities in the homogenized model are determined, standard numerical techniques can be used to compute the macroscopic response at a cost independent of the small scales. However, deriving such homogenized models often requires strong simplifying assumptions, e.g., periodicity of the underlying microstructure, which may not be seen too often in practical applications. This has led to the birth and development of more general purpose multiscale numerical methods over the last two decades, allowing for accurate approximations of the macroscale behaviour without fully resolving the microscale variations over much larger macroscopic geometries, see e.g., \cite{Abd09a,E11,HeM14,HFM98,MaP14}.

Realistic models may incorporate uncertainty (randomness) either because the model parameters are not perfectly known or because the intrinsic nature of the problem is uncertain. Moreover, the uncertainty may arise from multiple sources, e.g. from the coefficients, the right-hand side, and the initial data, or from the computational geometry itself. In this paper, we consider the following multiscale elliptic PDE with multiscale random coefficients, modelling e.g., diffusion or conductivity in random media over a domain $D\subset \R^d$:
\begin{equation}
\label{eq:model problem}
\left\lbrace
\begin{aligned}
-\nabla \cdot \left( a^{\varepsilon}(x) \nabla u^{\varepsilon} \right) &= f & \quad  & \text{in }D, \\
u^{\varepsilon} &= 0 & \quad & \text{on } \partial D ,
\end{aligned}
\right.
\end{equation}
where $a^{\varepsilon}(x) = a(x/\varepsilon)$ represents the heterogeneous medium with microscopic variations of characteristic length $\varepsilon \ll \abs{D} = \mathcal{O}(1)$. We assume that the tensor $a^{\varepsilon}(\cdot)$ is the realization of a random field $a^{\varepsilon}:\Omega\times D \mapsto \R^{d\times d}$, where $\Omega$ is a suitable probability space, and that $a^{\varepsilon}(\cdot) \in \left[ L^{\infty} \left(\Omega,L^{\infty}(D)\right)\right]^{d \times d}$ is symmetric, uniformly elliptic and bounded, i.e., $\exists \, \alpha,\, \beta > 0$ such that
\begin{equation}\label{eq: continuity coercivity assumption}
\alpha \abs{ \zeta }^2 \le \zeta \cdot a^{\varepsilon}(x) \zeta \le \beta \abs{ \zeta }^2 ,\, \forall \zeta \in \R^d ,\, \text{ a.e. in } \Omega\times D ,\, \forall \varepsilon>0.
\end{equation}
The original problem \cref{eq:model problem} is well-posed for any $f  \in H^{-1}(D)$. Under the condition that the coefficients are stationary and ergodic, the solution $u^{\varepsilon}$ of \cref{eq:model problem} can be approximated, as $\varepsilon \to 0$, by the solution $u^{0}$ of the so-called \emph{homogenized equation}:
\begin{equation} \label{eq:homogenized elliptic}
\left\lbrace
\begin{aligned}
-\nabla \cdot \left( a^{0} \nabla u^{0} \right) &= f & \quad & \text{in }D ,\\
u^{0} &= 0 & \quad & \text{on } \partial D ,
\end{aligned}
\right.
\end{equation}
where $a^0$ is a deterministic constant-valued matrix. The existence of the homogenized problem \eqref{eq:homogenized elliptic} is the central result of homogenization theory \cite{BLP78,CiD99,JKO94}. 

The first theoretical results on the homogenization of stationary ergodic random media were proved in \cite{Koz79,PaV79}. The homogenized tensor can be recovered by solving the \emph{corrector problem}:
\begin{equation}\label{eq: ell corr spat prob}
-\nabla\cdot\left(a(x)\left(\nabla\chi + \xi\right)\right) = 0 \quad\text{ in }\R^d,
\end{equation}
where $\xi\in \R^d, \abs{\xi}=1$. The gradient of the corrector function, $\nabla\chi$, is a stationary random field itself and the homogenized coefficient in \eqref{eq:homogenized elliptic} can be expressed as
\begin{equation}
\label{eq: a0}
\xi \cdot a^0\xi = \E{\left(\nabla \chi + \xi\right) \cdot a(x) \left(\nabla \chi + \xi\right)} = \fint_{\R^d} \left(\nabla \chi + \xi\right) \cdot a(x) \left(\nabla \chi + \xi\right) \,dx,
\end{equation}
where the $x$-dependency is lost due to the stationarity of the coefficients and the second identity is guaranteed by the Birkhoff's ergodic theorem. 
The auxiliary problem \eqref{eq: ell corr spat prob} is set over the entire $\mathbb{R}^d$, and cannot be computed numerically. Hence, the effective coefficients are approximated by solving approximate models, some of which are reviewed in \Cref{subsec: existing approaches}.

	\subsection{Existing elliptic models for computational homogenization}
	\label{subsec: existing approaches}
	In computations, a truncation of the infinite domain $\mathbb{R}^d$ in formulations \eqref{eq: ell corr spat prob}-\eqref{eq: a0} is inevitable. In what follows, we will mention three computational approaches to approximate the homogenized coefficient, two of which have already been analysed in the stationary-ergodic random setting. 
		\subsubsection{A naive elliptic approach}
		\label{subsec: naive elliptic}
		The most natural approach for approximating the homogenized coefficient is to replace the infinite domain formulations \eqref{eq: ell corr spat prob}-\eqref{eq: a0} by 
\begin{equation}\label{eq: ell corr spat prob domain R}
\left\lbrace
\begin{aligned}
-\nabla\cdot\left(a(x)\left(\nabla\chi_R + \xi\right)\right) &= 0 &\quad&\text{ in }K_R:=(-R/2,R/2)^d, \\
\chi_R &= 0 & \quad & \text{on } \partial K_R ,
\end{aligned}
\right.
\end{equation}
and 
\begin{equation}
\label{eq: a0RL}
\xi \cdot a^{0,R,L}\xi = \fint_{K_L} \left(\nabla \chi_R + \xi\right) \cdot a(x) \left(\nabla \chi_R + \xi\right) \,dx,
\end{equation}
with  $K_L:=(-L/2,L/2)^d\subset K_R$. In theory any $1\leq L \le R$ may be chosen in \eqref{eq: a0RL}, but for computational efficiency, $L = \mathcal{O}(R)$ is preferred in practice. We implicitly assume that $L=\mathcal{O}(R)$ although we keep the notation $L$ as a separate parameter to indicate that the size of the averaging domain $|K_L|$ may be different than $|K_{R}|$. Note that the Dirichlet boundary conditions (BCs) in \eqref{eq: ell corr spat prob domain R}, may be safely replaced by other BCs, e.g. periodic or homogeneous Neumann BCs. 
The approximation error between $a^{0,R,L}$ and $a^{0}$ is measured in a mean square sense as 
\begin{multline} \label{Ineq:Error_Split1}
\sup\limits_{\xi\in\R^d,\norm{\xi}=1}\E{\abs{\xi \cdot (a^{0,R,L}- a^0)\xi}) ^2}^{\frac{1}{2}}
\le \underbrace{\sup\limits_{\xi\in\R^d,\norm{\xi}=1}\E{\abs{\xi \cdot (a^{0,R,L}- a^{0,L})\xi} ^2}^{\frac{1}{2}} 
}_{\text{boundary error}}\\
+ \underbrace{\sup\limits_{\xi\in\R^d,\norm{\xi}=1}\sqrt{\var{\xi \cdot a^{0,L}\xi}} 
}_{\text{statistical error}} .
\end{multline}
Here, the term $a^{0,L}$ is a random variable approximating the homogenized coefficient, i.e. ${\E{a^{0,L}}=a^0}$, and it is computed by replacing $\chi_{R}$ in formula \eqref{eq: a0RL} by the solution $\chi$ to \eqref{eq: ell corr spat prob}. The first term in the right-hand side of \eqref{Ineq:Error_Split1} is due to the mismatch between the values of $\chi$ and $\chi_R$ on the boundary $\partial K_{R}$, while the second term accounts for the \emph{statistical} error, due to the averaging of a finite number of samples. The boundary error is known to be $\mathcal{O}(R^{-\alpha})$ with $0<\alpha\leq 1$, \cite{BoP04,Yur86}. The convergence rate is known to be of first order ($\alpha=1$) for periodic coefficients, while it is typically worse in the random case. In dimension $d=1$ the convergence rate is proven to be $\alpha=1/2$, while $\alpha\approx\frac{6}{25}$ in dimension $d=3$, \cite{EMZ05}. On the other hand, the statistical error scales as $\mathcal{O}(L^{-d/2})$, which is the optimal rate due to the central limit theorem. Therefore the boundary error dominates the  errors due to discretization and statistical averaging, which is also present in typical multiscale numerical methods, and better strategies are needed to bring this error down to practical tolerances of interests. 

Several interesting approaches have recently been proposed to reduce the above mentioned boundary error, see \cite{AAP19a,AAP19c,AAP19b,ArR16b,BlL10,CEL15b,Glo08,GlH16,HMS19}, where the core idea relies on modifying the corrector problem \eqref{eq: ell corr spat prob domain R} so that the artificial boundary conditions imposed on the boundary $\partial K_{R}$ has a minor influence on the accuracy of the approximate homogenized coefficient. In the remaining part of the introduction, we briefly mention two successful approaches from \cite{AAP19c,Glo08}, which use elliptic corrector problems.
		\subsubsection{An elliptic problem with zero-th order regularization}
		\label{subsec: zero-th elliptic}
		Another approach relies on the following regularized corrector problem:
\begin{equation}\label{eq: regul prob}
\begin{cases}
\frac{1}{T} \chi_{T,R} - \nabla\cdot\left(a(x)\left(\nabla\chi_{T,R} + \xi\right)\right) = 0 & \text{in } K_{R}, \\
\chi_{T,R} = 0 &  \text{on } \partial K_R .
\end{cases}
\end{equation}
This regularized problem was first introduced in \cite{PaV79} to prove the existence of $\nabla \chi$ satisfying \eqref{eq: ell corr spat prob} in the abstract space of stationary random fields. Later in \cite{Glo08}, it was used directly as a model for computational approximation of the homogenized tensor as follows
\begin{equation}
\label{eq: a0TL}
\xi \cdot a^{0,R,L,T}\xi = \fint_{K_L} \left(\nabla \chi_{T,R} + \xi\right) \cdot a(x) \left(\nabla \chi_{T,R} + \xi\right) \,dx.
\end{equation} 
In this case, the error between $a^{0,R,L,T}$ and the exact homogenized coefficient $a^{0}$ is split into three terms
\begin{multline} \label{Ineq:Error_Split2}
\sup\limits_{\xi\in\R^d,\norm{\xi}=1}\E{\abs{\xi \cdot (a^{0,R,L,T}- a^0)\xi}) ^2}^{\frac{1}{2}}
\le
\underbrace{\sup\limits_{\xi\in\R^d,\norm{\xi}=1}\E{\abs{\xi \cdot (a^{0,R,L,T}- a^{0,L,T})\xi} ^2}^{\frac{1}{2}} 
}_{\text{boundary error}}\\
+\underbrace{ \sup\limits_{\xi\in\R^d,\norm{\xi}=1} \sqrt{\var{\xi \cdot a^{0,T,L}\xi}}} _{\text{statistical error}} 
+\underbrace{\sup\limits_{\xi\in\R^d,\norm{\xi}=1} \abs{\xi \cdot(a^{0,T} - a^0) \xi}}_{\text{systematic error}},
\end{multline}
where the deterministic constant coefficient $a^{0,T}$ is defined in \eqref{eq: definition a0T Gloria} , and $a^{0,L,T}$ is obtained by replacing $\chi_{T,R}$ in \eqref{eq: a0TL} by $\chi_{T}$ solving \eqref{eq: regul prob} over $\mathbb{R}^{d}$, where the Dirichlet boundary conditions are (naturally) not imposed. 
These terms are denoted, respectively, by \emph{boundary}, \emph{statistical} and \emph{systematic} errors.
The first term in the right-hand side of \eqref{Ineq:Error_Split2} accounts for the \emph{boundary} error due to the mismatch between $\chi_{T,R}$ and $\chi_{T}$. 
On the other hand, the term $a^{0,T}$ is the expected value of $a^{0,L,T}$ that satisfies 
\begin{equation}
\label{eq: definition a0T Gloria}
\xi\cdot a^{0,T} \xi := \E{\xi\cdot a^{0,L,T} \xi} = \lim\limits_{L\to +\infty}\fint_{K_L} \left(\nabla \chi_{T} + \xi\right) \cdot a(x) \left(\nabla \chi_{T} + \xi\right) \,dx,
\end{equation}
by the Birkhoff's ergodic theorem.
Hence, the second term in the right-hand side of \eqref{Ineq:Error_Split2} is due to the \emph{statistical} error associated with the approximation \eqref{eq: a0TL}. 
Finally, the last term in the right hands of \eqref{Ineq:Error_Split2} represents the error due to adding a zero-th order regularization term of order $1/T$ to the standard elliptic corrector problem \eqref{eq: ell corr spat prob}, and is called the \emph{systematic} (or bias) error. 
The precise decay of the errors in relation with boundary, systematic, and statistical errors are summarized in Table \ref{Table_Comparison_Error_Estimates} to improve readability.
		\subsubsection{An elliptic problem with an exponential regularization}
		Recently, another approach based on a modified elliptic PDE was proposed in \cite{AAP19a} and analysed in \cite{AAP19c} for periodic coefficients. Similar to \eqref{eq: regul prob}, the approach relies on adding a regularization term to the standard elliptic corrector problem \eqref{eq: ell corr spat prob domain R} in order to reduce the effect of inaccurate homogeneous Dirichlet BCs posed on the boundary $\partial K_{R}$ of the computational geometry. In this case, the corrector problem reads as
 \begin{equation}\label{eq: modified_Elliptic}
\begin{cases}
e^{-\mathcal{A} T} g_{\xi} - \nabla\cdot\left(a(x)\left(\nabla\chi_{T,R} + \xi\right)\right) = 0 &\quad \text{in } K_{R}, \\
\chi_{T,R} = 0 & \quad  \text{on } \partial K_R ,
\end{cases}
\end{equation}
where $\mathcal{A} := -\nabla \cdot \left( a(x) \nabla \right)$, and $g_{\xi}:= \nabla \cdot \left(a(x) \xi\right)$. The approximate homogenized coefficient is then obtained by computing \eqref{eq: a0TL} as in the zero-th order regularization approach. The approximation error can again be split into three terms (boundary, systematic, and statistical errors) exactly as in \eqref{Ineq:Error_Split2}, whose decay rates are summarized in Table \ref{Table_Comparison_Error_Estimates}.
As it can be observed, this approach also leads to an exponential decay of the boundary error, see \cite{AAP19c} for a proof which is valid for periodic deterministic coefficients. Although both regularized approaches \eqref{eq: regul prob} and \eqref{eq: modified_Elliptic} lead to an exponential decay of the boundary error, the latter has advantages with respect to two aspects: 
\begin{enumerate}[label=\alph*)]
	\item in both approaches the optimal values for $L$ and $T$ are, respectively, $L = \mathcal{O}(R)$, and the $T = \mathcal{O}(R^\beta)$, $1<\beta<2$, but in approach \eqref{eq: regul prob} the exponential term decays as $e^{-c_1 R^{1-\alpha/2}}$, while it decays as $e^{-c_2 R^{2-\alpha}}$ for the approach \eqref{eq: modified_Elliptic};
	%
	%
	\item for dimensions up to $d=5$, the prefactor of the boundary error for the exponential regularization approach grows much slower than that of the zero-th order regularization approach. 
\end{enumerate}
Therefore, in simulations, the exponential decay of the boundary error for the model \eqref{eq: modified_Elliptic} is observed for moderate sizes of $R$, e.g., $R=10$, whereas in the zero-th order regularization approach \eqref{eq: regul prob} much larger values for $R$ are needed to bring the boundary error to computationally desirable tolerances. 
\renewcommand{\arraystretch}{1.5}
\begin{table}[h]
	\centering
	\begin{tabular}{| c | c  c  c |}
		\hline
		Corrector prob. & Naive \eqref{eq: ell corr spat prob domain R} &  zero-th order reg. \eqref{eq: regul prob}  & exponential reg. \eqref{eq: modified_Elliptic}  \\
		\hline
		Boundary err. & $R^{-\alpha}$   & $\sqrt{T} e^{- c_1 \frac{R-L}{\sqrt{T}}}$ & $\frac{R^{d-1} T^{(5-d)/2}}{|R-L|^3} e^{-c_2 \frac{|R-L|^2}{T}} $   \\
		Systematic err.  & $0$  & $T^{-d/2}$ & $T^{-d/2}$ (proved here)    \\
		Statistical err.   & $L^{-d/2}$  &  $L^{-d/2}$ & $L^{-d/2}$ \\
		\hline
	\end{tabular}
\caption {Error bounds for three different approaches. Note that the error estimates for the boundary error in relation with the exponential regularization approach is proved only in deterministic periodic setting.}
\label{Table_Comparison_Error_Estimates}
\end{table}
\renewcommand{\arraystretch}{1}
	\subsection{The goals of the paper}
	As discussed in the previous subsection, the corrector problem \eqref{eq: modified_Elliptic} has recently been found to be very practical from a ``cost vs accuracy'' point of view. However, the previous analysis in \cite{AAP19c} covered only periodically varying deterministic microstructures. The main goal of the present paper is to extend the part of the existing periodic theories to stationary ergodic random media, and thereby prove the validity of the approach for media possessing more complicated microstructure.  

In general, a complete analysis of a given corrector problem is linked to the following main theoretical questions:
\begin{itemize}
\item well-posedness of the corrector/equation,
\item study of the boundary error, 
\item study of the systematic error (if it exists),
\item study of the statistical error.
\end{itemize}
Several methods to improve the decay of the boundary error in the random homogenization setting are currently available. For example, the ``embedded'' method proposed in \cite{CEL15b,CEL20} approximates the homogenized coefficients by solving a variational minimization problem over a domain where a cell with heterogeneous coefficients is embedded into a homogeneous environment.
The study of the boundary error for the exponentially regularized corrector problem \eqref{eq: modified_Elliptic} is associated with the decay of the Green's function for the parabolic PDEs on bounded domains, which is reported in a previous paper of the authors, see \cite{AAP19c}. Although the previous analysis has been conducted for periodic tensors, in principle it is possible to prove exponentially decaying error estimates also when the periodicity assumption is relaxed. This can be achieved by following the proofs in the periodic setting together with additional regularity assumptions on the coefficient $a$. 
Therefore, the study of the boundary error is excluded from the present work. 
Moreover, for the statistical error, it is known that one can not perform better than the optimal decay rate $\mathcal{O}(L^{-d/2})$ originating from the central limit theorem, which is yet another theoretical consideration not included here. Several works addressed the issue of mitigating the statistical error, \cite{EKL15,LeM15b}. A control variate approach for reducing the variance of numerical approximations of the homogenized coefficients is been employed in \cite{LeM15b} in the context of periodic media with randomly positioned defects. In the case of non-ergodic media, an efficient approach to approximate $a^0$ relies on the Multi Level Monte Carlo method \cite{Gil08}, as proposed in \cite{EKL15}.

The goals of the present paper are to establish the well-posedness of the corrector problem \eqref{eq: modified_Elliptic} for stationary-ergodic tensors, and to study the systematic error due to adding the exponential regularization term. 
In particular, the second goal makes it possible to choose optimal values for the parameters $T,R,L$ in the model problem \eqref{eq: modified_Elliptic}, which is needed for computationally efficient and accurate approximations of the homogenized tensor in random media.

The paper is structured as follows. In \Cref{sec:notation} we describe our notation and recall some results that will be used in the exposition. 
We describe the modified elliptic approach in \Cref{sec:mod ell prob} and discuss and prove the main results for the stochastic setting in \Cref{sec:main results}. 
Numerical experiments are reported in \Cref{sec:numerical experiments}. 

\section{Notations and definitions}\label{sec:notation}	
In this section we explain the used notation and provide a precise formulation of the proposed corrector problems in the stochastic setting. 
\begin{itemize}
	\item Let $\Omega$ denote the set of measurable functions $a:\R^d\mapsto\R^{d\times d}$ such that:
	\begin{enumerate}[label=\roman*)]
		\item $a(x)$ is symmetric: $a_{ij}(x)=a_{ji}(x)$, and
		\item $a(\cdot)$ is uniformly elliptic and bounded, i.e., $\exists\alpha,\beta>0$ such that 
		\[
		\alpha\abs{\xi}^2 \le \xi \cdot a(x) \xi \le \beta \abs{\xi}^2,\quad \forall\xi\in\R^d,\text{ a.e. } x\in\R^d.
		\]
	\end{enumerate}
	\item Let $\mathcal{F}$ be the $\sigma$-algebra generated by the family
	\begin{equation*}
	\left\lbrace a\in\Omega\mapsto \int_{\R^d} a_{ij}(x)\varphi(x)\,dx, \varphi\in \Ccinf{\R^d},i,j=1,\dots,d\right\rbrace,
	\end{equation*}
	that means that two realizations $\check{a}$ and $\hat{a}$ can be identified in $\Omega$ if they differ only on a zero measure set.
	\item Let us endow the measurable space $\left(\Omega,\mathcal{F}\right)$ with a probability measure $\mathbb{P}$. The expected value of any random variable $X:\Omega\mapsto \R$ is denoted by $\mathbb{E}[X]:=\int_{\Omega} X\,d\mathbb{P}$.
\end{itemize}
\begin{definition}
	A \emph{translation group} (or $d$-dimensional dynamical system) is a family of invertible measurable maps, indexed by $x\in\R^d$, $\tau_x:\Omega \mapsto \Omega$ such that 
	\begin{enumerate}[label=\roman*)]
		\item $\tau_{x+y} = \tau_x\tau_y$, $\tau_0 = \text{Id}$;
		\item $\tau_x$ preserves the measure $\mathbb{P}$: $\mathbb{P}(\tau_x F) = \mathbb{P}(F)$, for any $F \in \mathcal{F}$ and any $x \in \R^d$;
		\item for any random variable $X:\Omega\mapsto \R$, the function $X(\tau_x a(\cdot))$ is $\R^d\times\Omega$ measurable with respect to the product $\sigma$-algebra.
	\end{enumerate}
\end{definition}
In the present work we will use the translation group $\tau_x$ defined by: 
\begin{equation*}
	\tau_x a(y) := a(x+y).
\end{equation*}
\begin{definition}
	Let $\mathcal{B}$ denote the Borel $\sigma$-algebra on $\R^d$. A \emph{stationary} random variable is a $\mathcal{B}\times\mathcal{F}$-measurable map $X: \R^d\times\Omega \mapsto \R$ such that, for any $y\in \R^d$
	\begin{equation*}
	X(x,\tau_y a(\cdot)) = X (x+y,a(\cdot)), \text{a.e.} x\in \R^d,\mathbb{P}\text{-a.e. } a\in \Omega.
	\end{equation*}
	If the random variable takes values in $\R^d$ it will be called \emph{random vector field}. If the random variable takes values in $\R^{d\times d}$ it will be called \emph{random tensor field}.
\end{definition}
\begin{remark}
	The identical random tensor field $X:\R^d\times\Omega \mapsto \R^{d\times d}$ that satisfies
	\begin{equation*}
	X(x,a(\cdot)) = a(x)
	\end{equation*}
	is stationary by definition. So, the random coefficients are stationary by definition. 
\end{remark}
\begin{definition}
	A set $F\in \mathcal{F}$ is called \emph{invariant} if $\tau_xF \subset F$ and a measurable function $X:\Omega\mapsto \R$ is called \emph{invariant} if $X(\tau_x \omega) = X(\omega)$ almost everywhere in $\Omega$. 
	
	A translation group $\tau_x$ is called \emph{ergodic} if the only invariant sets $F$ have either $\Pp(F)=0$ or $\Pp(F)=1$ or, alternatively, if all invariant functions are constant almost everywhere in $\Omega$.
\end{definition}
\begin{itemize}
	\item We denote by $K_R:=[-R/2,R/2]^d$ the hypercube with side length $R$.
	\item Let $U\subset \R^d$ be a domain. The Sobolev space $W^{p,k}(U)$ is defined as 
	$$W^{p,k}(U):= \{ f: D^{\gamma} f \in L^{p}(U) \text{ for all multi-index } \gamma \text{ with } |\gamma| \leq k\}.$$
	The norm of a function $f \in W^{p,k}(U)$ is given by 
	\begin{align*}
	\| f \|_{W^{p,k}(U)} := \begin{cases} \left( \sum_{|\gamma| \leq k} \int_{U} |D^{\gamma} f(x)|^{p} \; dx \right)^{1/p} & (1 \leq p < \infty) \\
	\sum_{|\gamma| \leq k} \text{ess sup}_{U} |D^{\gamma} f| & (p=\infty).
	\end{cases}
	\end{align*}
	\item We denote by $L^p_{loc}(U)$ the space of locally $p$-integrable functions,
	\begin{equation*}
	L^p_{loc}(U) = \left\lbrace f:U\mapsto \R: \int_K \abs{f(x)}^p\,dx<+\infty,\quad \forall K\subset U, K\text{ compact}\right\rbrace.
	\end{equation*}
	For any $f\in L^1_{loc}(\R^d)$ it is possible to define the quantity
	\begin{equation*}
	\fint_{\R^d} f(x)\,dx:= \lim\limits_{R\to +\infty} \fint_{K_R} f(x)\,dx = \E{f(x)},
	\end{equation*}
	where the identity follows from the ergodic theorem.
	\item We denote by $W^{p,k}_{loc}(U)$ the space of locally $W^{p,k}(U)$ functions,
	\begin{equation*}
	W^{p,k}_{loc}(U) = \left\lbrace f\in L^p_{loc}(U):  D^{\gamma} f \in L^{p}_{loc}(U) \text{ for all multi-index } \gamma \text{ with } |\gamma| \leq k \right\rbrace.
	\end{equation*} 
	The space $W^{2,1}_{loc}(U)$ will also be denoted by $H^1_{loc}(U)$.
	\item We denote by $L^p_{pot, loc}(U)$ the space of potential locally $p$-integrable potential functions,
	\begin{equation*}
	L^p_{pot, loc}(U) = \left\lbrace F\in \left(L^{p}_{loc}(U)\right)^d : \exists u \in W^{p,1}_{loc}(U), F = \nabla u \right\rbrace.
	\end{equation*}
	\item We denote by $\mathcal{L}^p$ the space of stationary fields,
	\begin{equation*}
	\mathcal{L}^p = \left\lbrace u \in L^p_{loc}(\R^d): u \text{ is stationary } \mathbb{P}\text{-a.s.}\right\rbrace.
	\end{equation*}
	The norm of a function $u \in \mathcal{L}^p$ is given by 
	\begin{equation*}
	\| u \|_{\mathcal{L}^p} := \left( \fint_{\R^d} \abs{u(x)}^p\,dx \right)^{\frac{1}{p}}.
	\end{equation*}
	\item We denote by $\mathcal{L}^p_{pot}$ the space of stationary potential fields,
	\begin{equation*}
	\mathcal{L}^p_{pot} = \left\lbrace F \in L^p_{pot, loc}(\R^d): F \text{ is stationary } \mathbb{P}\text{-a.s.}\right\rbrace.
	\end{equation*}
	The space $\mathcal{L}^p_{pot}$ can also be viewed as the completion in the  $\mathcal{L}^p$-norm of the set
	\begin{equation*}
	\left\lbrace \nabla u: u\in C^{\infty}(\R^d) \text{ is a stationary field} \right\rbrace.
	\end{equation*}
	It is remarkable that the function $u$ such that $\nabla u(\cdot,a) = F(\cdot,a)$ is not necessarily stationary. 
	\item We denote by $\mathcal{H}^1$ the space of stationary weakly differentiable fields,
	\begin{equation*}
	\mathcal{H}^1 = \left\lbrace u \in \mathcal{L}^2: \nabla u \in \mathcal{L}^2\right\rbrace.
	\end{equation*}
	The norm of a function $u \in \mathcal{H}^1$ is given by 
	\begin{equation*}
	\| u \|^2_{\mathcal{H}^1} := \fint_{\R^d} \abs{u(x)}^2\,dx + \fint_{\R^d} \abs{\nabla u(x)}^2\,dx.
	\end{equation*}
	For stationary random variables on ergodic environments one can equivalently write
	\begin{equation*}
	\| u \|^2_{\mathcal{H}^1} = \E{\abs{u(x)}^2} + \E{\abs{\nabla u(x)}^2}.
	\end{equation*}
	\item Let $I\subset R$ be an interval and $X$ be a normed space. We denote by $L^2(I,X)$ the space of square integrable functions $u:I\mapsto X$:
	\begin{equation*}
	L^2(I,X) = \left\lbrace u:I\mapsto X  : \int_{I} \norm{u(t)}_X^2\,dt<+\infty \right\rbrace
	\end{equation*}
	\item Let $I\subset R$ be an open interval and $X$ be a normed space. We denote by $C(I,X)$ the space of functions $u:I\mapsto X$ continuous for any $t\in I$.
	%
	%
\end{itemize}

A fundamental result of qualitative stochastic homogenization theory is the existence of the corrector functions $\chi$, as proved in \cite{PaV79}:
\begin{theorem}
	Let $a(\cdot)$ be a stationary and ergodic tensor field. Then, for any direction $\xi\in \R^d$, $\abs{\xi}=1$, there exists a unique $\chi \in L^2(\Omega,H^1_{loc}(\R^d))$ such that
	\begin{enumerate}[label=\roman*)]
		\item $\chi$ satisfies the corrector problem
		\begin{equation*}
		-\nabla\cdot\left(a(\nabla \chi + \xi)\right)=0, \text{ in }\mathcal{D}'(\R^d), \Pp\text{-a.s.};
		\end{equation*}
		\item $\displaystyle\chi(0,a) = 0$;
		\item $\displaystyle\E{\nabla \chi} = 0 $;
		\item $\nabla \chi$ is a stationary vector field.
	\end{enumerate}	
Moreover, $\chi$ grows sub-linearly at infinity: for every compact set $K\subset \R^d$,
\begin{equation*}
\lim\limits_{R\to\infty} \sup_{x\in K} \left(\frac{1}{R}\chi(Rx,a)\right)^2=0, \quad\forall a \in \Omega.
\end{equation*}
\end{theorem}
\begin{remark}
In the general case, the corrector $\chi$ is not statistically stationary. It is possible to prove (see, e.g. \cite[Chapter 4]{AKM19} or \cite{GlO17}) that the corrector $\chi$ is stationary for $d>2$.
\end{remark}

Throughout the exposition, we assume that $u$ is the solution of the following parabolic differential problem:
\begin{equation} \label{eq:parabolic problem}
\left\lbrace
\begin{aligned}
& \frac{\partial u}{\partial t} - \nabla\cdot(a(x)\nabla u) = 0 & \quad &\text{in } \R^d\times(0,+\infty),\\
& u(x,0) = \nabla \cdot (a(x)\xi) & \quad &\text{in } \R^d.
\end{aligned}
\right. 
\end{equation}
The solution of \eqref{eq:parabolic problem} is well defined by the use of the parabolic Green's function, see also \cite{AAP19a,AAP19c,AAP19b} and \cite[Chapter 9]{AKM19}:
\begin{equation}\label{eq: duhamel formula u}
u(x,t,a) = - \int_{\R^d} \nabla_y \Gamma (x,y,t) \cdot a(y)\xi \,dy,
\end{equation}
where $\Gamma (x,y,t) $ is the fundamental solution of the parabolic equation and solves:
\begin{equation*}
\left\lbrace
\begin{aligned}
& \frac{\partial \Gamma}{\partial t}(\cdot,y,\cdot) - \nabla_x\cdot(a(x)\nabla_x \Gamma(\cdot,y,\cdot)) = 0, \\
& \Gamma(\cdot,y,0) = \delta_y(\cdot), 
\end{aligned}
\right.
\quad
\text{and}
\quad
\left\lbrace
\begin{aligned}
& \frac{\partial \Gamma}{\partial t}(x,\cdot,\cdot) - \nabla_y\cdot(a(y)\nabla_y \Gamma(x,\cdot,\cdot)) = 0, \\
& \Gamma(x,\cdot,0) = \delta_x(\cdot), 
\end{aligned}
\right.
\end{equation*}
where $\delta_x$ is the \emph{Dirac delta function} centered in $x$. From classical results on linear parabolic partial differential equations \cite{Fri64}, we know that, for any $a\in \Omega$, problem \eqref{eq:parabolic problem} has a unique solution $u$.

\section{The modified elliptic model}\label{sec:mod ell prob}
In the context of deterministic homogenization, the use of corrector problem \eqref{eq: ell corr spat prob} on a reference volume element $K_R$ revealed to be computationally inefficient for approximating the true homogenized coefficients $a^0$. Indeed, for growing values of $R$, the number of unknowns grows as $R^d$, while the accuracy of the approximation of $a^0$ scales as $R^{-1}$, \cite{BoP04,Yur86}. 
In \cite{AAP19a,AAP19c} a higher order numerical homogenization scheme was presented, analysed and numerically tested for the case of periodic coefficients.
The novel upscaling method is based on \emph{modified} elliptic cell problems and in this paper we will study its accuracy properties in the framework of homogenization of stochastic coefficients. Let us thus consider the cell problem
\begin{equation}\label{eq: spatial mod ell model KR}
\left\lbrace
\begin{aligned}
&-\nabla \cdot \left(a(x) \left(\nabla \chi_{R,T} + \xi \right)\right) = - e^{-\mathcal{A}T}\left[\nabla\cdot\left(a(x)\xi\right)\right]\,, \quad \text{in }K_R,\\
&\chi_{R,T}\text{ is $K_R$-periodic,}
\end{aligned}
\right.
\end{equation}
where $\mathcal{A}$ is the second order elliptic operator endowed with periodic boundary conditions defined by $\mathcal{A}:=-\nabla\cdot\left(a(x)\nabla\right)$ and $e^{-\mathcal{A}T}$ is the semigroup generated by $\mathcal{A}$ evaluated at time $T$. 
As for the regularized cell problem \eqref{eq: regul prob}, we believe that the modified corrector over the bounded cell $\chi_{R,T}$ approximates with an infinite order of accuracy the modified corrector over the unbounded cell $\chi_T$, that is defined as the solution of 
\begin{equation}\label{eq: spatial mod ell model}
-\nabla \cdot \left(a(x) \left(\nabla \chi_{T} + \xi \right)\right) = - u(\cdot,T)\,, \quad \text{in }\R^d,
\end{equation}
where $u(\cdot,T)$ is the solution of the Cauchy problem \eqref{eq:parabolic problem} evaluated at time $T$. 
In this perspective, the approximation of $\chi_{T}$ over a bounded cell lies within a negligible error away from $\chi_{T}$, provided that the cell size $R^d$ is sufficiently large. This result has been proved in \cite{AAP19c} for regular, periodic coefficients.  

The modified corrector functions $\chi_{R,T}$ and $\chi_T$ are employed to upscale the multiscale tensor by the \emph{cell average formulas}:
\begin{gather}
\label{eq: definition a0RLT}
\xi\cdot a^{0,R,L,T} \xi: = \fint_{K_L} \left(\nabla \chi_{R,T} + \xi\right)\cdot a(x) \left(\nabla \chi_{R,T} + \xi\right) \,dx,\text{ and }\\
\label{eq: definition a0LT}
\xi\cdot a^{0,L,T} \xi: = \fint_{K_L} \left(\nabla \chi_{T} + \xi\right)\cdot a(x) \left(\nabla \chi_{T} + \xi\right) \,dx,
\end{gather}
where the restriction over the smaller box $K_L$ is necessary in order to achieve the exponential decay of the boundary error, whose upper bound depends on $(R-L)$, as proved in \cite{AAP19c} in the periodic setting.
The main contribution of this work is to prove the well-posedness of the modified cell problem on the unbounded domain \eqref{eq: spatial mod ell model}, prove a convergence rate for the systematic error $a^{0,T}-a^0$ ($a^{0,T}$ being defined in \eqref{eq: identity a0T}) and demonstrate the decay of the global modelling error $a^{0,R,L,T}-a^0$ by means of numerical experiments.
%


\section{Main result}\label{sec:main results}
In this section we provide the main results of the present work, formulated in \Cref{thm: well posedness2,thm: systematic error}. 
As stated in \Cref{sec:mod ell prob}, our goal is to bound the systematic errors for the upscaling scheme \eqref{eq: spatial mod ell model}. First of all, we prove that the cell problem \eqref{eq: spatial mod ell model} is well-posed and that $\nabla \chi_T$ is a stationary random field. 
The well-posedness proof is based on the equivalence between the gradient of the modified corrector $\nabla\chi_T$ and the time integral of $\nabla u$, for which we rely on time decay properties of parabolic solutions.
The stationarity of $\nabla \chi_T$, which can be compared to the stationarity of $\nabla \chi$ of \eqref{eq: ell corr spat prob}, is essential for applying the ergodic theorem in the definition of $a^{0,T}$. Then, we can exploit the time decay of $\nabla u$ as well as the equivalent definition of $a^{0,T}$ and $a^{0}$ as integrals in $\R^d$ and $\Omega$ to derive \textit{a-priori} bounds on the systematic error.

\begin{theorem}\label{thm: well posedness2}
	Let $u\in C\left([0,+\infty),\mathcal{L}^2\right)$ be the solution of \eqref{eq:parabolic problem}. 
	Then, there exists a unique $\nabla \chi_T\in \mathcal{L}^2_{pot}$ such that
	\begin{equation}\label{eq: mod ell corr prob}
	-\nabla\cdot\left(a(\nabla \chi_T + \xi)\right)= - u(\cdot,T), \text{ in }\mathcal{D}'(\R^d),\, \mathbb{P}\text{-a.s.}.
	\end{equation}
\end{theorem}
\begin{remark}
	In \cite{AKM19} the authors proved that for $\nabla u$ and $\nabla\chi$ as above the following relation holds true:
	\begin{equation*}
	\nabla \chi(x) = \int_{0}^{+\infty} \nabla u(x,t)\,dt,\quad \mathbb{P}\text{-a.s.}
	\end{equation*}
	This identity entails the stationarity of $\nabla \chi$, as a consequence of the stationarity of $\nabla u$. Moreover, the identity does not hold true for $u$ and $\chi$, because the time integral of $u$ does not converge.
	We will use a similar identity to prove \Cref{thm: well posedness2}. 
\end{remark} 
%
%
%
%
\begin{theorem}[systematic error]\label{thm: systematic error}
	Let $a(x)\in \Omega$, $a^{0}$ and $a^{0,T}$ be defined, respectively, as in \eqref{eq: a0} and by
	\begin{equation} \label{eq: identity a0T}
	\xi \cdot a^{0,T} \xi 
	:=\lim\limits_{L\to + \infty}\fint_{K_L} \left(\nabla \chi_T + \xi \right) \cdot a(y) \left(\nabla \chi_T + \xi \right) \,dy
	= \E{\xi\cdot a^{0,T,L} \xi}.
	\end{equation}
	Then, there exists a positive constant $C(d, \alpha,\beta)<+\infty$ such that
	\begin{equation}\label{eq: thesis thm systematic}
	\sup_{\xi\in\R^d,\norm{\xi}=1}\abs{\xi \cdot (a^{0,T} - a^0) \xi } \le C T^{-\frac{d}{2}}.
	\end{equation}
\end{theorem}
The proofs of \Cref{thm: well posedness2,thm: systematic error} are based on the decay in time of the parabolic solution $u$, which are collected in \Cref{sec:parabolic}.
	\subsection{Decay of parabolic solutions}\label{sec:parabolic}
	In this section we collect some results about the decay in time of the solutions to parabolic PDEs in $\R^d$ with stationary random coefficients. 

First of all, we recall a classical result on the time decay of the solutions to parabolic problems and deduce the results of Lemma \ref{lemma: stationarity u}. These results are not new, for example they are proved in \cite{AKM19} for the case of $\mathbb{Z}^d$-stationary random fields. 
\begin{lemma}\label{lemma:classical result time decay}
	Let $u$ be the solution of \eqref{eq:parabolic problem}. Then, there exists a constant $C(\alpha,\beta,d)>0$ such that, for every $t>0$,
	\begin{equation}
	\norm{u(\cdot,t)}_{L^{\infty}(\R^d)} + t^{\frac{1}{2}} \norm{\nabla u(\cdot,t)}_{L^{\infty}(\R^d)} \le C t^{-\frac{1}{2}}.
	\end{equation}
\end{lemma} 
\begin{lemma}\label{lemma: stationarity u}
	Let $u$ be the solution of \eqref{eq:parabolic problem}. Then, $u$ is a stationary random field and 
	\begin{equation}
	\E{u(x,t)} = 0, \quad \forall t>0, \forall x\in \R^d.
	\end{equation}
\end{lemma}
\begin{proof}
	\textit{Step 1.} We prove the stationarity of $u$. Let us recall that $u$ can be expressed by formula \eqref{eq: duhamel formula u}. Then, by the fact that $\Gamma (x+z,y+z,t,a)=\Gamma (x,y,t,\tau_z a)$
	\begin{align*}
	u(x+z,t,a) &= - \int_{\R^d} \nabla_y \Gamma (x+z,y,t,a) \cdot a(y)\xi \,dy\\
	&= - \int_{\R^d} \nabla_y \Gamma (x+z,y+z,t,a) \cdot a(y+z)\xi \,dy\\
	&= - \int_{\R^d} \nabla_y \Gamma (x,y,t,\tau_z a) \cdot \tau_z a(y)\xi \,dy = u(x,t,\tau_z a).
	\end{align*}
	\textit{Step 2.} Let $B_1$ be the unit ball centered at $0$, $\psi \in \Ccinf{B_1}$ with unit mass in $L^1(B_1)$ and $\psi_R(x):= R^{-d}\psi(x/R)$. Let us write \eqref{eq:parabolic problem} in weak form with $\psi_R$ as test function and integrate in time for $0<t_1<t<t_2$:
	\begin{equation*}
	\E{u(\cdot,t_1)} - \E{u(\cdot,t_2)} = \lim\limits_{R \to +\infty}\E{\int_{t_1}^{t_2} \int_{\R^d} \nabla u(x,t) \cdot a(x) \nabla \psi_R(x) \,dx\,dt}. 
	\end{equation*}
	By the H\"older inequality, we bound the absolute value of the right-hand side from above:
	\begin{equation*}
	\begin{split}
	\abs{\E{\int_{t_1}^{t_2} \int_{\R^d} \nabla u(x,t) \cdot a(x) \nabla \psi_R(x) \,dx\,dt}} \le
	\E{\int_{t_1}^{t_2} \beta \norm{\nabla u(\cdot,t)}_{L^{\infty}(\R^d)} \norm{\nabla \psi_R}_{L^1(\R^d)} \,dt}\\
	\le \beta R^{-1} \norm{\nabla \psi}_{L^1(\R^d)} \E{\int_{t_1}^{t_2} \norm{\nabla u(\cdot,t)}_{L^{\infty}(\R^d)}  \,dt}
	\end{split}
	\end{equation*}
	The term $\E{\int_{t_1}^{t_2} \norm{\nabla u(\cdot,t)}_{L^{\infty}(\R^d)}  \,dt}$ is uniformly bounded in $R$ thanks to the decay of $\norm{\nabla u(\cdot,t)}_{L^{\infty}(\R^d)}$ of Lemma \ref{lemma:classical result time decay}. So,
	\begin{equation*}
	\lim\limits_{R \to +\infty}\E{\int_{t_1}^{t_2} \int_{\R^d} \nabla u(x,t) \cdot a(x) \nabla \psi_R(x) \,dx\,dt} = 0
	\end{equation*}
	and we deduce that $\E{u(\cdot,t)}$ is constant in time. From the fact that $\norm{u(\cdot,t)}_{L^{\infty}(\R^d)}$ decays to zero and from the stationarity of $u$, we conclude that $\E{u(x,t)}=0$ for any $t>0$ and any $x\in\R^d$ . 
\end{proof}

Time decay rates of $\E{\abs{u}^{p}}$ and $\E{\abs{\nabla u}^{p}}$ for homogenization problems over discrete networks were proved in several works, e.g. \cite[Theorem 1]{GNO15} and \cite[Lemma 9.7]{Mou19}:
\begin{equation*}
\E{\abs{u}^{p}}^{\frac{1}{p}} \le C (t+1)^{-\left( \frac{1}{2} + \frac{d}{4}\right)}\text{ for any } p\ge1
\quad \text{ and }\quad
\E{\abs{\nabla u}^{2}}^{\frac{1}{2}} \le C (t+1)^{-\left( 1 + \frac{d}{4}\right)}.
\end{equation*}

More recently, similar estimates were also derived for the continuous case in \cite{AKM19}. 
\Cref{thm: time decay u} and Corollary \ref{cor: boundedness E u} provide time decay bounds on the moments $\E{\abs{u}^p}$. 
\begin{theorem}[{\cite[Theorem 9.1]{AKM19}}]\label{thm: time decay u}
	For every $\sigma\in (0,2)$, there exists a constant $C(\sigma,d,\alpha, \beta)<+\infty$ such that the following holds. Let $a(\cdot)\in L^{\infty}(\R^d)$ be a stationary random field such that, for every $x\in\R^d$, $a(x)\xi$ is $\mathcal{F}$-measurable and let $u\in\mathcal{P}$ be the solution of \eqref{eq:parabolic problem}. Then, for every $t\in[1,+\infty)$ and $x\in\R^d$, 
	\begin{equation}\label{eq: Mourrat time decay u}
	\E{\exp\left( \left(C^{-1} t^{\frac{1}{2} + \frac{d}{4}} \abs{u(x,t)}\right)^{\sigma} \right)} \le 2.
	\end{equation}
\end{theorem}
\begin{corollary}\label{cor: boundedness E u}
	Let the assumptions of Theorem \ref{thm: time decay u} be satisfied. Then, for any $p\ge 1$, there exists a constant $C(p,d,\alpha,\beta)<+\infty$ such that, for every  $t\in[1,+\infty)$ and $x\in\R^d$
	\begin{equation}\label{eq:decay u p}
	\E{\abs{u(x,t)}^p}^\frac{1}{p} \le C t^{-\left(\frac{1}{2} + \frac{d}{4}\right)}
	\end{equation}
\end{corollary}
\begin{proof}
	From Theorem \ref{thm: time decay u}, by taking $\sigma=1$, we know that there exists a constant $C(d,\alpha, \beta)<+\infty$ such that 
	\begin{equation*}
	\E{\exp\left( C^{-1} t^{\frac{1}{2} + \frac{d}{4}} \abs{u(x,t)} \right)} \le 2
	\end{equation*}
	for every $t\in[1,+\infty)$ and $x\in\R^d$. Since the exponential of a random variable $X$ grows faster than $\abs{X}^p$ for any $p$, the integrability of $e^X$ implies the integrability of any power of $X$. Therefore, there exists a constant $C(p)<+\infty$ such that 
	\begin{equation*}
	\E{\abs{X}^p} \le C(p) \E{e^X}. 
	\end{equation*}
	By taking $X= C^{-1} t^{\frac{1}{2} + \frac{d}{4}} \abs{u(x,t)}$ in the previous inequality we conclude that there exists a constant $C(p,d,\alpha,\beta)<+\infty$ such that 
	\begin{equation*}
	\E{\abs{u(x,t)}^p}^\frac{1}{p} \le C t^{-\left(\frac{1}{2} + \frac{d}{4}\right)}.
	\end{equation*}
	\hfill
\end{proof}

Corollary \ref{cor: boundedness E u} shows that there is a clear difference between the time decay of parabolic solutions set in bounded domains (as, for instance, in the case of periodic correctors) and in unbounded domains (as in the stochastic homogenization setting).
Indeed, in the periodic (or bounded domain) setting, the Poincar\'e inequality entails exponential decay in time of the spatial $L^2$-norm. Such a property is fundamental in the derivation of exponential order convergence rates of the modelling error in \cite{AAP19a,AAP19c,AAP19b}. In the stochastic setting we do not necessarily have such an inequality in $\mathcal{H}^1$.

\begin{proposition}\label{prop: u is L2H1 C0L2}
Let $u$ be the solution of \eqref{eq:parabolic problem} with $\nabla\cdot\left(a(x)\xi\right)\in \mathcal{L}^2$. Then
\begin{equation*}
u\in L^2\left((0,+\infty),\mathcal{H}^1\right) \cap C\left([0,+\infty),\mathcal{L}^2\right).
\end{equation*}
\end{proposition}
\begin{proof}
	We first prove that $u\in L^2\left((0,+\infty),\mathcal{H}^1\right)$ and, then, that $u\in C\left([0,+\infty),\mathcal{L}^2\right)$.\\
	\textit{Step 1} - $u\in L^2\left((0,+\infty),\mathcal{H}^1\right)$:
	
	We already know from \Cref{lemma: stationarity u} that $u(\cdot,t)$ is stationary for any $t\ge 0$. So, we only have to prove that 
	\begin{equation*}
		\int_{0}^{+\infty}\E{u(\cdot,t)^2}\,dt<+ \infty \text{ and } \int_{0}^{+\infty}\E{\abs{\nabla u(\cdot,t)}^2}\,dt<+ \infty.
	\end{equation*}
	The function $\E{u(\cdot,t)^2}$ is decreasing in time, indeed, from \eqref{eq:parabolic problem},
	\begin{equation*}
	\frac{d}{dt}  \E{u^2} = 2  \E{u\partial_t u} = -2 \E{\nabla u \cdot a(x) \nabla u} <0.
	\end{equation*}	 
	So, we can bound the integral using $\nabla\cdot\left(a(x)\xi\right)\in \mathcal{L}^2$ and the result of Corollary \ref{cor: boundedness E u}:
	\begin{equation}\label{eq: u in L2}
	\begin{aligned}
	\int_{0}^{+\infty} \E{u(\cdot,t)^2} \,dt &\le \int_{0}^{1} \E{u(\cdot,t)^2} \,dt + \int_{1}^{+\infty} \E{u(\cdot,t)^2} \,dt\\
	&\le \E{\abs{\nabla\cdot\left(a(x)\xi\right)}^2} + C \int_{1}^{+\infty} t^{-\left(1 + \frac{d}{2}\right)} \,dt < + \infty.
	\end{aligned}
	\end{equation}
	Next, from the ellipticity of $a(\cdot)$, we have:
	\begin{equation*}
	\alpha \E{\abs{\nabla u(\cdot,t)}^2} \le \E{\nabla u \cdot a(x) \nabla u} = -\frac{1}{2} \frac{d}{dt}  \E{u(\cdot,t)^2}. 
	\end{equation*}
	So, since $\E{u(\cdot,t)^2}$ vanishes for $t \to + \infty$,
	\begin{equation}\label{eq: grad u in L2}
	\int_{0}^{+\infty} \E{\abs{\nabla u(\cdot,t)}^2} \,dt \le \frac{1}{2\alpha} \E{u(\cdot,0)^2}  < + \infty.
	\end{equation}
	From \eqref{eq: u in L2} and \eqref{eq: grad u in L2} we conclude that $u\in L^2\left((0,+\infty),\mathcal{H}^1\right)$.
	
	\textit{Step 2} -  $u\in C\left([0,+\infty),\mathcal{L}^2\right)$:
	
	Let $t\ge0$. Since $f(z)=\sqrt{z}$ is continuous in $[0,+\infty)$, it is sufficient to prove the continuity of $\E{u^2}$: 
	\begin{equation*}
	\begin{aligned}
	 \E{u(\cdot,t+h)^2} - \E{u(\cdot,t)^2} &= \int_{t}^{t+h} \frac{d}{dt}\E{u(\cdot,t)^2} \,dt \\
	 &= -\int_{t}^{t+h} \E{\nabla u \cdot a(x)\nabla u} \,dt \underset{h\to 0}{\longrightarrow} 0,
	\end{aligned}
	\end{equation*}
	and the proof is concluded.
\end{proof}


Now, we state a result on the time decay of the second moment of $\nabla u$. The proof follows from the one of \cite[Lemma 9.7]{Mou19}.
\begin{proposition}\label{prop: time decay DU 2}
	Let $a(\cdot)\in \Omega$ and let $u$ be the solution of \eqref{eq:parabolic problem}. Then, there exists a positive constant $C(d,\alpha,\beta) < + \infty$ such that, for every $t\in[2,+\infty)$ and $x\in\R^d$,
	\begin{equation}
	\label{eq: time decay Du 2}
	\E{\abs{\nabla u(x,t)}^{2}}^{\frac{1}{2}} \le C t^{-(\frac{d}{4}+1)} .
	\end{equation}
\end{proposition}
\begin{proof}
	Let us begin by proving that the map $t\mapsto\E{\nabla u(x,t)\cdot a(x) \nabla u(x,t)}$ is nonincreasing. Indeed, its time derivative can be expressed as:
	\begin{align*}
	\partial_t \E{\nabla u(x,t)\cdot a(x) \nabla u(x,t)} &= 2 \E{ \nabla (\partial_t u)(x,t)\cdot a(x) \nabla u(x,t)}\\ 
	&= 2 \E{ \nabla (\nabla \cdot \left(a(x) \nabla u(x,t)\right))\cdot a(x) \nabla u(x,t)}\\ 
	&= -2 \E{ \abs{\nabla \cdot \left(a(x) \nabla u(x,t)\right)}^2}\le 0.
	\end{align*}  
	Thus, from the weak formulation of \eqref{eq:parabolic problem} with $u$ as test function and inequality \eqref{eq:decay u p} for $t/2\ge 1$ and $p=2$, we can write
	\begin{align*}
	\E{\nabla u(x,t) \cdot a(x) \nabla u(x,t)} &\le \frac{2}{t} \int_{\frac{t}{2}}^{t} \E{\nabla u(x,s) \cdot a(x) \nabla u(x,s) } \,ds\\
	&\le -\frac{1}{t} \int_{\frac{t}{2}}^{t} \partial_t \E{\abs{u(x,s)}^2} \,ds\\
	&\le \frac{1}{t} \left[\E{\abs{u\left(x,t\right)}^{2}} - \E{\abs{u\left(x,\frac{t}{2}\right)}^{2}}\right] \\
	&\le C t^{-(\frac{d}{2}+2)}.
	\end{align*}
	Then, \eqref{eq: time decay Du 2} follows from the assumption of uniform ellipticity of the coefficients.
\end{proof}

	\subsection{Well-posedness modified corrector problem}\label{sec:well posedness}
Now we are ready to prove that problem \eqref{eq: spatial mod ell model} is well-posed.
\begin{proof}[\Cref{thm: well posedness2}]
	Let us define the stationary function
	\begin{equation}\label{eq: Psi def}
	\Psi := \int_{0}^{T} \nabla u(\cdot,t)\,dt.
	\end{equation}
	Then, $\Psi\in (\mathcal{L}^2)^d$. Indeed, by Minkowski integral inequality and Proposition \ref{prop: u is L2H1 C0L2} we have:
	\begin{equation*}
	\E{\abs{\Psi(x)}^2}^{\frac{1}{2}} 
	:= \E{\abs{\int_{0}^{T} \nabla u(x,t)\,dt}^2}^{\frac{1}{2}}
	\le  \int_{0}^{T} \E{\abs{\nabla u(x,t)}^2}^{\frac{1}{2}}\,dt < + \infty.
	\end{equation*}
	The weak form of \eqref{eq:parabolic problem} is: Find $u\in L^2\left((0,+\infty),\mathcal{H}^1\right)$ such that 
	\begin{equation*}
	\frac{d}{dt} \E{u\phi} + \E{\nabla \phi \cdot a(x) \nabla u} = 0 ,\quad \forall \phi \in \mathcal{H}^1. 
	\end{equation*}
	By integration in time and \eqref{eq: Psi def}, we get
	\begin{equation*}
	\E{\nabla \phi \cdot a(x) \left(\Psi + \xi\right) }   = -\E{u(\cdot,T)\phi} ,\quad \forall \phi \in \mathcal{H}^1.
	\end{equation*}
	To conclude, we have to prove that $\Psi \in \mathcal{L}^2_{pot}$.
	The function $\Psi$ is trivially vortex-free, since it is the gradient of $\int_{0}^{T} u(\cdot,t)\,dt$. Hence, we are allowed to define $\Psi$ as $\nabla \chi_T$.
	
	The uniqueness of $\nabla \chi_T$ trivially follows from uniqueness of solution for the standard corrector problem \eqref{eq: ell corr spat prob}, proved in \cite{PaV79}.
\end{proof}
	\subsection{Proof of the systematic error bound}\label{sec:bias error proof}
	The systematic error is 
\begin{equation}
e_{SYS} := \sup_{\xi\in\R^d,\norm{\xi}=1} \abs{\xi \cdot (a^{0,T} - a^0) \xi }.
\end{equation}
We will rely on the result on the time decay of $\nabla u(\cdot,t)$ and on the definition of $\nabla \chi$ and $\nabla \chi_T$ as time integral of $\nabla u(\cdot,t)$ in order to bound the systematic error.
\begin{proof}[Theorem \ref{thm: systematic error}]
	We first notice that the two identities contained in \eqref{eq: identity a0T} follow from the stationarity of $a(\cdot)$ and the Birkhoff ergodic theorem. 
	Next, we prove that 
	\begin{equation}
	\xi\cdot (a^{0,T} - a^0) \xi
	= \E{ \left(\nabla \chi_T - \nabla \chi \right)\cdot a(\cdot) \left(\nabla \chi_T - \nabla \chi \right)}.
	\end{equation}
	By definition of $a^{0,T}$ and $a^0$,
	\begin{equation*}
	\begin{aligned}
	\xi\cdot (a^{0,T} - a^0) \xi
	&= \E{ \left(\nabla \chi_T +\xi \right)\cdot a(\cdot) \left(\nabla \chi_T +\xi \right) - \left(\nabla \chi +\xi \right)\cdot a(\cdot) \left(\nabla \chi +\xi \right)}\\
	&= \E{ \left(\nabla \chi_T +\xi \right)\cdot a(\cdot) \left(\nabla \chi_T +\xi \right)  - \left(\nabla \chi +\xi \right)\cdot a(\cdot) \left(\nabla \chi_T +\xi \right)}\\ 
	&\quad + \E{ \left(\nabla \chi_T +\xi \right)\cdot a(\cdot) \left(\nabla \chi +\xi \right) - \left(\nabla \chi +\xi \right)\cdot a(\cdot) \left(\nabla \chi +\xi \right)}\\
	&= \E{ \left(\nabla \chi_T - \nabla \chi \right)\cdot a(\cdot) \left(\nabla \chi_T +\xi \right) + \left(\nabla \chi_T - \nabla \chi \right)\cdot a(\cdot) \left(\nabla \chi +\xi \right) }\\
	&= \E{ \left(\nabla \chi_T - \nabla \chi \right)\cdot a(\cdot) \left(\nabla \chi_T +\xi \right) - \left(\nabla \chi_T - \nabla \chi \right)\cdot a(\cdot) \left(\nabla \chi +\xi \right) }\\
	&= \E{ \left(\nabla \chi_T - \nabla \chi \right)\cdot a(\cdot) \left(\nabla \chi_T -\nabla \chi \right) },
	\end{aligned}
	\end{equation*}
	where the fourth inequality comes from 
	\begin{equation*}
	\E{\left(\nabla \chi_T - \nabla \chi \right)\cdot a(x) \left(\nabla \chi +\xi \right) } = 0 = -\E{\left(\nabla \chi_T - \nabla \chi \right)\cdot a(x) \left(\nabla \chi +\xi \right) },
	\end{equation*}
	for any $x\in \R^d$\footnote{To prove it, one can test the standard corrector equation against $\theta(\chi_T - \chi)$, where $\theta\in \Ccinf{\R^d}$ is such that $\theta\equiv 1$ on the ball $B_R$, then, pass to the limit for $R\to +\infty$ and use the ergodic theorem.}. 
	Thus, by the uniform boundedness of $a(\cdot)$ and the H\"older inequality we have:
	\begin{equation*}
	\abs{\xi\cdot \left(a^{0,T} - a^0\right)\xi} 
	= \abs{\E{ \left(\nabla \chi_T-\nabla \chi\right) a(\cdot) \left(\nabla \chi_T-\nabla \chi\right)}}
	\le \beta \E{\abs{\nabla \chi_T-\nabla \chi}^2}.
	\end{equation*}
	Now, we recall that 
	\begin{equation*}
	\nabla \chi = \int_{0}^{+\infty} \nabla u(\cdot,t) \,dt, \quad \text{ and } \nabla \chi_T = \int_{0}^{T} \nabla u(\cdot,t) \,dt,
	\end{equation*}
	we subsitute these equivalences in the expression above and use the Minkowski integral inequality to switch the two integrations: 
	\begin{equation*}
	\abs{\xi\cdot \left(a^{0,T} - a^0\right)\xi} 
	\le \beta \E{\abs{\int_{T}^{+\infty} \nabla u(\cdot,t) \,dt}^2}
	\le \beta \left(\int_{T}^{+\infty}  \E{\abs{\nabla u(\cdot,t)}^2}^{\frac{1}{2}} \,dt\right)^2.
	\end{equation*}
	Finally, from the time decay result for $\nabla u(\cdot,t)$ \eqref{eq: time decay Du 2} we conclude that 
	\begin{equation*}
	\abs{\xi\cdot \left(a^{0,T} - a^0\right)\xi} 
	\le C \left(\int_{T}^{+\infty} t^{-\left(\frac{d}{4}+1\right)} \,dt\right)^2
	\le C T^{-\frac{d}{2}}.
	\end{equation*}
	and \eqref{eq: thesis thm systematic} follows from the fact that $C$ does not depend on $\xi$. 
\end{proof}

\section{Numerical experiments} \label{sec:numerical experiments}
In this section we collect the results of numerical experiments performed in order to:
\begin{itemize}
	\item verify numerically the exponential convergence of the boundary error in the context of random coefficients;
	\item verify numerically the correctness of the proved bound on the systematic error in Theorem \ref{thm: systematic error};
	\item compare the convergence of the global modelling error for the standard numerical homogenization scheme and for the modified elliptic approach.
\end{itemize}
The convergence of the error is measured with respect to the modelling parameters $R$, $L$ and $T$. 
As we already saw in \Cref{sec: intro}, the mean square error can be bounded by the sum of three terms: the boundary, the systematic and the statistical errors.
In order to avoid that the approximation error is dominated by the statistical error (which is often the case when $L$ is not chosen large enough), we looked at the error in mean, i.e.,
\begin{equation}\label{eq:mean err}
e_{M} := \sup\limits_{\xi\in\R^d,\norm{\xi}=1} \abs{  \xi \cdot \E{a^{0,R,L,T}- a^0}\xi  }.
\end{equation}
The advantage of using this error measure is that it does not depend on the statistical error, allowing for experimentally testing the decay of the sum of boundary and systematic errors. This is clear from the decomposition: 
\begin{equation*}
\abs{ \xi \cdot \E{a^{0,R,L,T}- a^0}\xi } 
\le 
\abs{ \xi \cdot \E{a^{0,R,L,T}- a^{0,L,T}}\xi } 
+ \underbrace{ \abs{ \xi \cdot \E{a^{0,L,T}- a^{0,T}}\xi } }_{=0\text{, by \eqref{eq: identity a0T} } }
+ \abs{ \xi \cdot \E{a^{0,T}- a^0}\xi },
\end{equation*}
thus implying that 
\begin{equation*}
e_{M} \le \underbrace{  \sup\limits_{\xi\in\R^d,\norm{\xi}=1} \abs{ \xi \cdot \E{a^{0,R,L,T}- a^{0,L,T}}\xi } }_{\text{boundary error} }
+ \underbrace{  \sup\limits_{\xi\in\R^d,\norm{\xi}=1} \abs{ \xi \cdot \left(a^{0,T}- a^{0}\right)\xi } }_{\text{systematic error} }.
\end{equation*}
The boundary error has a deterministic upper bound that converges exponentially with exponent $\frac{-\abs{R-L}^2}{T}$, see \Cref{Table_Comparison_Error_Estimates}, while the systematic error decays as $T^{-d/2}$, as proved in \Cref{sec:main results}. Computing $e_M$ \emph{exactly} is not possible, because it relies on an integral formulation over the probability space $\Omega$, but we approximate it by computing the empirical average
\begin{equation}
\label{eq: empirical average}
\bar{a}^{0,R,L,T,N} = \frac{1}{N} \sum_{k=1}^{N} a^{0,R,L,T,k},
\end{equation}
converging to $\E{a^{0,R,L,T}}$ in probability, by the weak law of large numbers. Upon choosing a sufficiently large number of samples, $N$, the difference between $\bar{a}^{0,R,L,T,N} $ and $\E{a^{0,R,L,T}} $ can be made negligible. 
	\subsection{The covariance function of random fields}
\begin{figure}[h!]
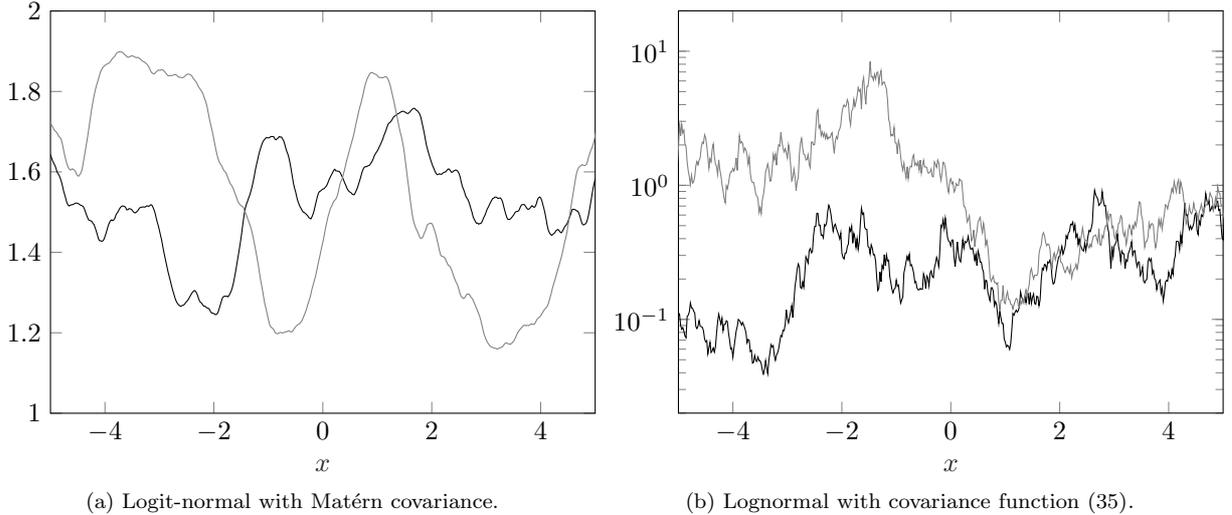

	\subfloat[Logit-normal with Mat\'ern covariance.]{
		\input{\imgdir/logitNormal1D/coeff_Matern}}
	\subfloat[Lognormal with covariance function \eqref{eq: invSqrt cov func}.]{
		\input{\imgdir/lognormal1D/coeff_invSqrt}}
	\caption{Realizations of random fields.}
	\label{fig:random fields}
\end{figure}

	For simplicity, we will consider only isotropic media, for which the heterogeneous tensor can be written as $a(x) = f(x)I$, where $f:\R^d\mapsto \R$ and $I\in \R^{d\times d}$ is the identity matrix. With a slight abuse of notation we will denote both the matrix of coefficients and the function $f$ above by $a(x)$.
Stationarity of the random fields implies that $\E{a(\cdot)}=\mu I_{d\times d}$ does not depend on the spatial variable $x$ and that the covariance (matrix-valued) function $\mathrm{Cov}(x,y)$ defined as 
\begin{equation*}
\mathrm{Cov}(x,y) := \E{a(x)a(y)} - \mu^2
\end{equation*}
only depends on the distance $x-y$. Therefore, for stationary random fields, there exists a function $r:\R^d\mapsto \R$ such that 
\begin{equation*}
\mathrm{Cov}(x,y) = r(x-y).
\end{equation*}
When $r(t) = r(\abs{t})$ the medium is said to be \emph{statistically isotropic}.
Several choices for the covariance function are possible. For example, widely used classes of covariance functions for one dimensional isotropic random fields are the exponential covariance function:
\begin{equation*}
r(t) = \sigma^2 e^{-\abs{\frac{t}{l}}},
\end{equation*}
and the Mat\`ern covariance function:
\begin{equation*}
r(t) = \sigma^2 \frac{1}{\Gamma(\nu)2^{\nu-1}} \left(\sqrt{2\nu} \abs{\frac{t}{l}}\right)^{\nu} K_{\nu}\left(\sqrt{2\nu} \abs{\frac{t}{l}}\right),
\end{equation*}
where $\sigma^2$ is the variance, $l$ is the correlation length, $\Gamma$ is the gamma function, $K_{\nu}$ is the modified Bessel function of the second kind and $\nu$ is a smoothness parameter. 
Another choice is the the long-range covariance function:
\begin{equation} 
\label{eq: invSqrt cov func}
r(t) = \left( 1 + \abs{t}\right)^{-1/2}.
\end{equation}
All random fields considered in the numerical experiments are generated by the circulant embedding method described in \cite{DiN97}.
Two types of random fields are depicted in \Cref{fig:random fields}.
	\subsection{Optimal scaling of $T$ vs. $R$}
	Here we briefly discuss the optimal scaling of the parameter $T$ as a function of $R$ (and $L$) with the aim of maximizing the rate of decay of $e_M$.
We find the regime under which none of the boundary and systematic errors is dominating but the two are (approximately) equal. 
Let us start by considering the exponential term in the upper bound for the boundary error. Then, we impose:
\begin{equation*}
C_1 \exp\left(-c_2\frac{\abs{R-L}^2}{T}\right) = C_2 T^{-\frac{d}{2}} \implies T \log\left(\left(\frac{C_1}{C_2}\right)^{\frac{2}{d}} T\right) = \frac{2 c_2}{d}\abs{R-L}^2.
\end{equation*}
The constants $C_1,C_2$ are unknown and problem-dependent, but we can conclude that the optimal scaling is obtained for 
\begin{equation}
\label{eq: quasi-optimal T scaling}
c\abs{R-L}\le T \le C\abs{R-L}^2,
\end{equation}
with $c,C>0$.
	\subsection{One dimensional logit-normal random coefficients}
	We test the convergence of the approximate homogenized coefficient for a random diffusion coefficient distributed according to the logit-normal law. A logit-normal random field is an isotropic random field $a(\cdot)\in \Omega$ of the form 
\begin{equation*}
a(x) =\frac{b+e^{-\kappa(Z(x)-z_0)}}{c+e^{-\kappa(Z(x)-z_0)}};
\end{equation*}
where $b,c,\kappa,z_0\in\R$ and $Z$ is a Gaussian random field of zero average and $r(\cdot)$ is the covariance function. We set $b=2$, $c=1$, $\kappa=1$ and $x_0 = 0$. 
In the one dimensional case, the homogenized coefficient can be computed by the harmonic mean: $a^0 = \E{a(\cdot)^{-1}}^{-1}$. Hence, in the logit-normal case, 
\begin{equation*}
a^0 = \left( \int_{\R} a(y)^{-1} f_Z(y) \,dy\right)^{-1},
\end{equation*}
where $f_Z$ is the Gaussian probability density function.

We computed the approximation to the homogenized coefficients by finite elements (FE) discretization on a grid with mesh size $h=2^{-8}$. The modified auxiliary problem \eqref{eq: spatial mod ell model KR} is solved over the domain $K_R:=(-R/2,R/2)$ with periodic boundary conditions, with the values of $R$ ranging from 5 to 500. The other parameters are ${L=2R/3}$, for the size of the averaging domain $K_L$, and $T$, for the modified forcing term. 
The approximation to the homogenized coefficients are computed as in \eqref{eq: definition a0RLT}.
As an approximation of the quasi-optimal scaling \eqref{eq: quasi-optimal T scaling}, we choose 
$$ T = \frac{\abs{R-L}^2}{100}.$$


The right-hand side of \eqref{eq: spatial mod ell model KR} is approximated in the FE space by the exponential matrix $e^{-M_h^{-1}AT}\mathbf{g}$, where $\mathbf{g}$ is the vector of components of the projection of $g(\cdot) = \frac{d}{dx}a(\cdot)$ in the FE space, $M_h$ is the lumped mass matrix and $A_h$ is the stiffness matrix. The exponential matrix is not computed exactly, but it is approximated in the Krylov subspace generated by $M_h^{-1}A_h\mathbf{g}$ and computed by the Lanczos method ($M_h^{-1}A_h$ is symmetric and positive definite) as proposed in \cite{HoL97}. The maximum number of Krylov basis elements is 2000. 

\begin{figure}[h!]
	\centering
%
%
\definecolor{mycolor1}{rgb}{0.00000,0.44700,0.74100}%
\definecolor{mycolor2}{rgb}{0.85000,0.32500,0.09800}%
\definecolor{mycolor3}{rgb}{0.46667,0.67451,0.18824}%
\begin{tikzpicture}

\begin{axis}[%
width=.85\linewidth,
height=2.1in,
scale only axis,
xmode=log,
xmin=3,
xmax=600,
xminorticks=true,
xlabel style={font=\color{white!15!black}},
xlabel={$R$},
ymode=log,
ymin=5e-5,
ymax=.1,
yminorticks=true,
axis background/.style={fill=white},
ylabel={$\abs{ \bar{a}^{0,R,L,T,N} - a^{0}}$},
xmajorgrids,
ymajorgrids,
legend style={font=\scriptsize, legend cell align=left, align=left, draw=white!15!black,at={(1,1)},anchor=north east }
]
\addplot [color=mycolor1, mark=*, mark size=1pt, mark options={solid, mycolor1}]
  table[row sep=crcr]{%
5	0.0433211367787743\\
10	0.0283358838691772\\
15	0.0201301486675771\\
20	0.0141779683790313\\
25	0.0183196071900107\\
30	0.018982037852104\\
35	0.0166186568605025\\
40	0.0100875964044829\\
45	0.0078691528345185\\
50	0.00829037185110582\\
55	0.00563979873391807\\
60	0.00768321405043015\\
65	0.00576140769624955\\
70	0.00659977456099692\\
75	0.00463064391196677\\
80	0.0049034080350725\\
85	0.0031234561892195\\
90	0.00377444027641727\\
95	0.00112667564335744\\
100	0.00235762965525232\\
120	0.000870846817704551\\
140	0.000135873143604304\\
160	0.000906359911233112\\
180	0.000693565494072956\\
200	0.00264498991125772\\
220	0.000227578464816025\\
240	0.000870571260492836\\
260	0.000189095844625342\\
280	0.00105390740624678\\
300	6.85782165850135e-05\\
320	0.000218736356826454\\
340	9.59985520527784e-05\\
360	0.000426698796776703\\
380	0.000569236642575799\\
400	0.000617568099048116\\
420	0.000156652218083586\\
440	0.000590956856950209\\
460	0.000152435131200113\\
480	0.000264070485468393\\
500	0.000298970391377607\\
};
\addlegendentry{$L=\frac{2R}{3}, T = \frac{\abs{R-L}^2}{100}$}

\addplot [color=mycolor1, dashed]
  table[row sep=crcr]{%
5	0.06\\
10	0.03\\
15	0.02\\
20	0.015\\
25	0.012\\
30	0.01\\
35	0.00857142857142857\\
40	0.0075\\
45	0.00666666666666667\\
50	0.006\\
55	0.00545454545454545\\
60	0.005\\
65	0.00461538461538462\\
70	0.00428571428571429\\
75	0.004\\
80	0.00375\\
85	0.00352941176470588\\
90	0.00333333333333333\\
95	0.00315789473684211\\
100	0.003\\
120	0.0025\\
140	0.00214285714285714\\
160	0.001875\\
180	0.00166666666666667\\
200	0.0015\\
220	0.00136363636363636\\
240	0.00125\\
260	0.00115384615384615\\
280	0.00107142857142857\\
300	0.001\\
320	0.0009375\\
340	0.000882352941176471\\
360	0.000833333333333333\\
380	0.000789473684210526\\
400	0.00075\\
420	0.000714285714285714\\
440	0.000681818181818182\\
460	0.000652173913043478\\
480	0.000625\\
500	0.0006\\
};
\addlegendentry{$\mathcal{O}(R^{-1})$}

\addplot [color=mycolor2, mark=x, mark options={solid, mycolor2}]
  table[row sep=crcr]{%
5	0.0217584337863652\\
10	0.0146873819486577\\
15	0.000684064308650978\\
20	0.000220160421152693\\
25	0.00245990291932729\\
30	0.00410294933913824\\
35	0.00202653096664918\\
40	0.00173067468163479\\
45	0.00365027182180611\\
50	0.00118252423931353\\
55	0.00312108733349858\\
60	0.00280192241603383\\
65	0.00194916750786267\\
70	0.00160146548276496\\
75	0.000633699434974266\\
80	0.00340678766172542\\
85	8.84625107646642e-05\\
90	0.00316660789070689\\
95	0.00121801032731739\\
100	0.00119091129971238\\
120	0.000436303267045579\\
140	0.00103967577391839\\
160	0.00071385405374591\\
180	0.000640987035364748\\
200	0.000551891023467865\\
220	0.0012539632612556\\
240	0.000223286077615459\\
260	0.00117125462736745\\
280	0.00243688998881009\\
300	0.000266951105826641\\
320	0.000349768249598093\\
340	0.000758135002482785\\
360	0.00104524005119244\\
380	0.000164469971782433\\
400	0.00136279438360698\\
};
\addlegendentry{$L=R, T = \infty$ (standard)}

\addplot [color=mycolor2, dashed]
  table[row sep=crcr]{%
5	0.0134164078649987\\
10	0.00948683298050514\\
15	0.00774596669241483\\
20	0.00670820393249937\\
25	0.006\\
30	0.00547722557505166\\
35	0.0050709255283711\\
40	0.00474341649025257\\
45	0.00447213595499958\\
50	0.00424264068711928\\
55	0.00404519917477945\\
60	0.00387298334620742\\
65	0.00372104203767625\\
70	0.00358568582800318\\
75	0.00346410161513775\\
80	0.00335410196624968\\
85	0.00325395686727984\\
90	0.00316227766016838\\
95	0.00307793505625546\\
100	0.003\\
120	0.00273861278752583\\
140	0.00253546276418555\\
160	0.00237170824512628\\
180	0.00223606797749979\\
200	0.00212132034355964\\
220	0.00202259958738973\\
240	0.00193649167310371\\
260	0.00186052101883813\\
280	0.00179284291400159\\
300	0.00173205080756888\\
320	0.00167705098312484\\
340	0.00162697843363992\\
360	0.00158113883008419\\
380	0.00153896752812773\\
400	0.0015\\
420	0.0014638501094228\\
440	0.00143019388386839\\
460	0.00139875721236047\\
480	0.00136930639376292\\
500	0.00134164078649987\\
};
\addlegendentry{$\mathcal{O}(R^{-1/2})$}

\addplot [color=black, mark size=2.0pt, mark=triangle*, mark options={solid, black}]
table[row sep=crcr]{%
	10	0.0243770148782478\\
	20	0.0165035496808958\\
	50	0.00871694142088586\\
	100	0.00206026887286415\\
	500	0.000189568068751145\\
};
\addlegendentry{$L=\frac{2R}{3}, T = \frac{\abs{R-L}^2}{100}, N=10^5$}

\end{axis}
\end{tikzpicture}%
	\caption{Comparison of the error in mean of the standard and modified elliptic methods for a logit-normal random field with Mat\`ern covariance function of order $\nu = 3/2$.}
	\label{fig:mean_err_logit_matern}
\end{figure}
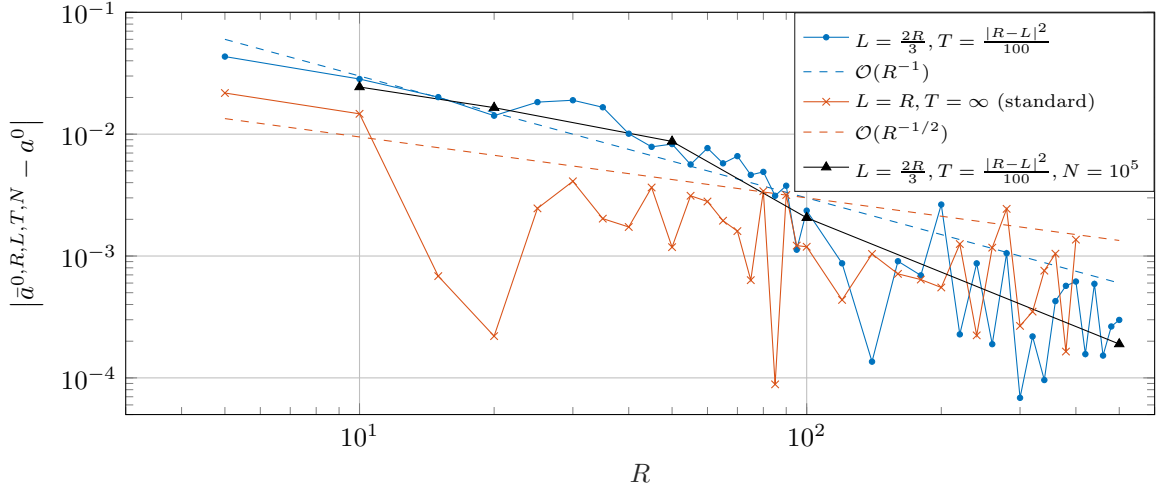

The error in mean between the approximate and the exact homogenized coefficient \eqref{eq:mean err} is plotted in \Cref{fig:mean_err_logit_matern}. Since the expected value of $a^{0,R,L,T}$ cannot be computed exactly, we approximate it by the empirical average \eqref{eq: empirical average} with $N=1000$ i.i.d. samples.
The red line shows the error decay for the standard auxiliary problem with periodic BCs and no oversampling. In this case, the only source of error is due to the BCs. The error for the modified elliptic approach is represented by the blue line. In this other case, the global error is the contribution of the boundary and systematic errors. 
As one can see, the red curve in \Cref{fig:mean_err_logit_matern} decays at a slow rate of $R^{-\frac{1}{2}}$.
The blue curve in \Cref{fig:mean_err_logit_matern} has a faster decay  of order (approximately) $\mathcal{O}(R^{-1})$, thanks to the scaling of $T$.
\Cref{fig:mean_err_logit_matern} also displays the error in mean when $N=10^5$ Monte Carlo samples are chosen. In these simulations, we make sure that the errors do not depend on the number of samples.

	\subsection{One dimensional lognormal random coefficients}
	We test the convergence of the approximate homogenized coefficient for a random diffusion coefficient distributed according to the lognormal law. A lognormal random field is an isotropic random field $a(\cdot)\in \Omega$ for which there exist $b,c>0$ such that
\begin{equation*}
a(x) = c e^{bZ(x)}, 
\end{equation*}
where $Z$ is a Gaussian random field of zero average and $r(\cdot)$ is the covariance function.
The model of lognormally distributed random coefficients is widely used in the environmental engineering community, see e.g. \cite{Dur91,EAW04,McC95}. 
However, such a coefficient does not belong to $\mathcal{M}(\alpha,\beta)$, so it is not guaranteed that it follows the theoretical estimates that we derived in the previous sections.
In the one dimensional case, the homogenized coefficient can be computed by the harmonic mean: $a^0 = \E{a(\cdot)^{-1}}^{-1}$. Hence, in the lognormal case, 
\begin{equation*}
a^0 = \left( \int_{\R} a(y)^{-1} f_Z(y) \,dy\right)^{-1} = c e^{-b^2/2},
\end{equation*}
where $f_Z$ is the Gaussian probability density function.
For this test, we have chosen ${b=c=1}$. 

We computed the approximation to the homogenized coefficients by finite elements (FE) discretization on a grid with mesh size $h=2^{-8}$. The modified auxiliary problem \eqref{eq: spatial mod ell model KR} is solved over the domain is $K_R:=(-R/2,R/2)$ with periodic boundary conditions. 
The values of the parameter $R$ vary from $5$ to $500$. The other parameters are ${L=2R/3}$, for the size of the averaging domain $K_L$, and $T$, which is chosen as 
\begin{equation*}
T\log(T) = \frac{\abs{R-L}^2}{100}.
\end{equation*} 
Note that since the lognormal random field does not belong to the class $\mathcal{M}(\alpha,\beta)$, it is not possible to choose the value of $T$ according to the optimal scaling derived in \eqref{eq: quasi-optimal T scaling}. The approximation to the homogenized coefficients are computed as in \eqref{eq: definition a0RLT}. Moreover, the right-hand side of \eqref{eq: spatial mod ell model KR} is approximated in the FE space as in the previous example.

\begin{figure}[h!]
	\centering
	\subfloat[$| \bar{a}^{0,R,L,T,N} - a^{0}|$.\label{fig:mean_err_lognorm}]{
%
%
\definecolor{mycolor1}{rgb}{0.00000,0.44700,0.74100}%
\definecolor{mycolor2}{rgb}{0.85000,0.32500,0.09800}%
\begin{tikzpicture}

\begin{axis}[%
width=.45\linewidth,
height=2.1in,
scale only axis,
xmode=log,
xmin=3,
xmax=600,
xminorticks=true,
xlabel style={font=\color{white!15!black}},
xlabel={$R$},
ymode=log,
ymin=0.02,
ymax=1,
yminorticks=true,
axis background/.style={fill=white},
xmajorgrids,
ymajorgrids,
legend style={font=\scriptsize, legend cell align=left, align=left, draw=white!15!black,at={(0.02,0.02)},anchor=south west }
]

\addplot [color=mycolor1, mark=*, mark size=1pt, mark options={solid, mycolor1}]
  table[row sep=crcr]{%
5	0.854432244637786\\
10	0.873408818265216\\
15	0.746727826446273\\
20	0.642718713120399\\
25	0.706264704960023\\
30	0.542295744408176\\
35	0.533670159437101\\
40	0.450432505543463\\
45	0.44838467594345\\
50	0.414891204397355\\
55	0.377936100267321\\
60	0.367468833612183\\
65	0.308315278539038\\
70	0.299003739773468\\
75	0.275901064471231\\
80	0.309029206965361\\
85	0.255543784245649\\
90	0.263443214706925\\
95	0.214233858179977\\
100	0.225177823735427\\
120	0.168527600977082\\
140	0.175611800873291\\
160	0.154609152170094\\
180	0.126868256381191\\
200	0.128745344014872\\
220	0.116390172078292\\
240	0.153152776924602\\
260	0.107373620718118\\
280	0.112241460273701\\
300	0.0890619319954424\\
320	0.0963243717625942\\
340	0.106687874235301\\
360	0.0842451381747734\\
380	0.106138402384822\\
400	0.101073610163996\\
420	0.0920919148265608\\
440	0.0975598952899854\\
460	0.0870587416791517\\
480	0.0802164614178291\\
500	0.0815706832619096\\
};
\addlegendentry{$L=\frac{2R}{3}$, $T\log(T) = \frac{\abs{R-L}^2}{100}$}

\addplot [color=mycolor1, dashed]
  table[row sep=crcr]{%
5	2.99069756244244\\
10	1.77827941003892\\
15	1.31199311417695\\
20	1.05737126344056\\
25	0.894427190999916\\
30	0.780115773106905\\
35	0.694942651170796\\
40	0.628716714841468\\
45	0.575560014246967\\
50	0.531829589694499\\
55	0.495139966423909\\
60	0.46385961395229\\
65	0.436832540682104\\
70	0.413215372645583\\
75	0.392377460851028\\
80	0.373837195305305\\
85	0.357220032554213\\
90	0.34223003202678\\
95	0.328630067510752\\
100	0.316227766016838\\
120	0.275812576637239\\
140	0.245699330589314\\
160	0.222284926254865\\
180	0.203491194526928\\
200	0.18803015465432\\
220	0.175058413947413\\
240	0.163999139272119\\
260	0.154443625879632\\
280	0.146093696044109\\
300	0.138726381676261\\
320	0.132171407930071\\
340	0.126296353697382\\
360	0.120996588185913\\
380	0.116188274619323\\
400	0.111803398874989\\
420	0.10778616799614\\
440	0.104090355703677\\
460	0.100677315362991\\
480	0.0975144716383632\\
500	0.0945741609003176\\
};
\addlegendentry{$\mathcal{O}(R^{-\frac{3}{4}})$}

\addplot [color=mycolor2, mark=x, mark options={solid, mycolor2}]
  table[row sep=crcr]{%
5	0.656952773361357\\
10	0.440700756140684\\
15	0.316118025492357\\
20	0.363461319049516\\
25	0.337870551093803\\
30	0.279467582564987\\
35	0.36498546678218\\
40	0.240346409900689\\
45	0.238495053056838\\
50	0.197310590020424\\
55	0.23992573733102\\
60	0.215951747627814\\
65	0.299670405846412\\
70	0.242508202625575\\
75	0.278312978188215\\
80	0.210180054340422\\
85	0.234426924818911\\
90	0.193825392249429\\
95	0.149754923521672\\
100	0.160257386520254\\
120	0.187762334109152\\
140	0.17437405211358\\
160	0.16706938614855\\
180	0.115468537747196\\
200	0.146745482034461\\
220	0.148814859049019\\
240	0.10682706555573\\
260	0.0935734641951442\\
280	0.106413570889513\\
300	0.116034413721469\\
320	0.12169503160519\\
340	0.115359089438903\\
360	0.0855565252396203\\
380	0.076399124313535\\
400	0.0860184077538798\\
420	0.0885775878811088\\
440	0.0876361159466698\\
460	0.0860338087295552\\
480	0.0658950415838799\\
500	0.101214926247539\\
};
\addlegendentry{$L=R, T =\infty$ (standard)}

\addplot [color=mycolor2, dashed]
  table[row sep=crcr]{%
  5	0.894427190999916\\
  10	0.632455532033676\\
  15	0.516397779494322\\
  20	0.447213595499958\\
  25	0.4\\
  30	0.365148371670111\\
  35	0.338061701891407\\
  40	0.316227766016838\\
  45	0.298142396999972\\
  50	0.282842712474619\\
  55	0.269679944985297\\
  60	0.258198889747161\\
  65	0.248069469178417\\
  70	0.239045721866879\\
  75	0.23094010767585\\
  80	0.223606797749979\\
  85	0.216930457818656\\
  90	0.210818510677892\\
  95	0.205195670417031\\
  100	0.2\\
  120	0.182574185835055\\
  140	0.169030850945703\\
  160	0.158113883008419\\
  180	0.149071198499986\\
  200	0.14142135623731\\
  220	0.134839972492648\\
  240	0.129099444873581\\
  260	0.124034734589208\\
  280	0.119522860933439\\
  300	0.115470053837925\\
  320	0.111803398874989\\
  340	0.108465228909328\\
  360	0.105409255338946\\
  380	0.102597835208515\\
  400	0.1\\
  420	0.0975900072948533\\
  440	0.0953462589245592\\
  460	0.0932504808240314\\
  480	0.0912870929175277\\
  500	0.0894427190999916\\
};
\addlegendentry{$\mathcal{O}(R^{-\frac{1}{2}})$}


\end{axis}
\end{tikzpicture}%
}
	\subfloat[$| \bar{a}^{0,R,L,T,N} - \bar{a}^{0,R_{max},L,T,N}|$.\label{fig:bd_err_lognorm}]{
%
%
\definecolor{mycolor1}{rgb}{0.00000,0.44700,0.74100}%
\definecolor{mycolor2}{rgb}{0.85000,0.32500,0.09800}%
\begin{tikzpicture}

\begin{axis}[%
width=.45\linewidth,
height=2.1in,
scale only axis,
xmode=log,
xmin=3,
xmax=600,
xminorticks=true,
xlabel style={font=\color{white!15!black}},
xlabel={$R$},
ymode=log,
ymin=0.0003,
ymax=3,
yminorticks=true,
axis background/.style={fill=white},
xmajorgrids,
ymajorgrids,
legend style={font=\scriptsize, legend cell align=left, align=left, draw=white!15!black,at={(0.02,0.02)},anchor=south west }
]

\addplot [color=mycolor1, mark=*, mark size=1pt, mark options={solid, mycolor1}]
table[row sep=crcr]{%
	5	1.12842533695664\\
	10	1.41245631613222\\
	15	1.50289753685928\\
	20	1.44301649169377\\
	25	1.47038876253139\\
	30	1.02709881478774\\
	35	0.892035411568606\\
	40	0.760561340796873\\
	45	0.663341060967759\\
	50	0.570109495898439\\
	55	0.523520048078749\\
	60	0.392973205436824\\
	65	0.344316371229227\\
	70	0.336044245253027\\
	75	0.341054688960483\\
	80	0.320431590784706\\
	85	0.234794730972734\\
	90	0.260882868002061\\
	95	0.179056487820767\\
	100	0.171388287840584\\
	120	0.0581251736446772\\
	140	0.0794570105197541\\
	160	0.00240792843537591\\
	180	0.0626709225254353\\
	200	0.0285741679675801\\
	220	0.0232422322983165\\
	260	0.028080271947136\\
	280	0.0114609533466644\\
	300	0.0364082540110644\\
	320	0.0240801390119194\\
	340	0.00986649526564781\\
	360	0.00977466897040119\\
	380	0.0138783500192191\\
	400	0.000450425799649518\\
	420	0.00155325362944544\\
	440	0.0200627363287068\\
	460	0.00707946399169745\\
	480	0.00692295211509833\\
	500	0.00261648302977791\\
};
\addlegendentry{$L=\frac{2R}{3}$, $T\log(T) = \frac{\abs{R-L}^2}{100}$}

\addplot [color=mycolor2, mark=x, mark options={solid, mycolor2}]
  table[row sep=crcr]{%
5	0.59654808784442\\
10	0.333782710161893\\
15	0.204097426883134\\
20	0.289376921763525\\
25	0.234476612429366\\
30	0.195648666684074\\
35	0.283535227929116\\
40	0.182559493667405\\
45	0.140046407721343\\
50	0.105602818637833\\
55	0.160252626967883\\
60	0.157762587006953\\
65	0.222262654846842\\
70	0.181022653347005\\
75	0.182011316366744\\
80	0.140420420377644\\
85	0.168457099759664\\
90	0.128652561870105\\
95	0.0609091189742381\\
100	0.0746589776147535\\
120	0.0854601401060467\\
140	0.115229362152346\\
160	0.0711939783841034\\
180	0.0367469354595177\\
200	0.0694491622468983\\
220	0.0659817949217753\\
240	0.00485883576522816\\
260	0.016206481105398\\
280	0.0298203702263076\\
300	0.0292948141991737\\
320	0.0230219269626443\\
340	0.0370287597779871\\
360	0.0019861731037476\\
380	0.00899004926458724\\
400	0.00399708731690651\\
420	0.01126545376149\\
440	0.0262641540299285\\
460	0.0145120840571954\\
480	0.00945215308834946\\
500	0.0337912696981053\\
};
\addlegendentry{$L=R$, $T=\infty$ (standard)}

\addplot [color=mycolor2, dashed]
  table[row sep=crcr]{%
5	2\\
10	1\\
15	0.666666666666667\\
20	0.5\\
25	0.4\\
30	0.333333333333333\\
35	0.285714285714286\\
40	0.25\\
45	0.222222222222222\\
50	0.2\\
55	0.181818181818182\\
60	0.166666666666667\\
65	0.153846153846154\\
70	0.142857142857143\\
75	0.133333333333333\\
80	0.125\\
85	0.117647058823529\\
90	0.111111111111111\\
95	0.105263157894737\\
100	0.1\\
120	0.0833333333333333\\
140	0.0714285714285714\\
160	0.0625\\
180	0.0555555555555556\\
200	0.05\\
220	0.0454545454545455\\
240	0.0416666666666667\\
260	0.0384615384615385\\
280	0.0357142857142857\\
300	0.0333333333333333\\
320	0.03125\\
340	0.0294117647058824\\
360	0.0277777777777778\\
380	0.0263157894736842\\
400	0.025\\
420	0.0238095238095238\\
440	0.0227272727272727\\
460	0.0217391304347826\\
480	0.0208333333333333\\
500	0.02\\
};
\addlegendentry{$\mathcal{O}(R^{-1})$}

\end{axis}
\end{tikzpicture}%
}
\caption{Error terms for a lognormal random field with covariance function $r(t) = \left(1 + \abs{t}\right)^{-1/2}$.}
\label{fig:err_lognorm}
\end{figure}
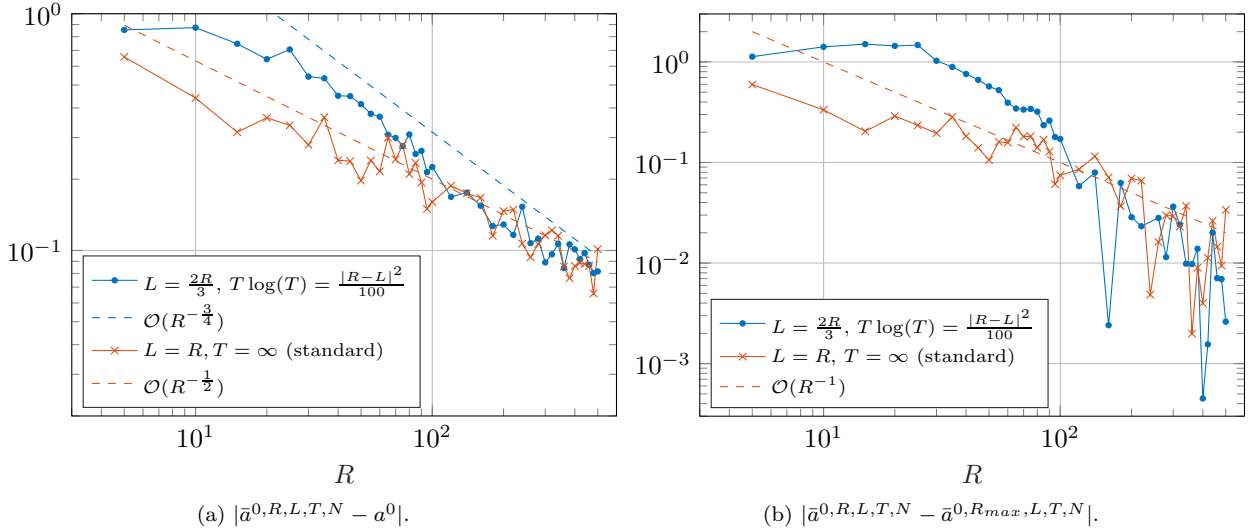

The error in mean between the approximate and the exact homogenized coefficient \eqref{eq:mean err} is plotted in \Cref{fig:mean_err_lognorm}. The expected value of $a^{0,R,L,T}$ is approximated it by  $N=1000$ independent samples of the lognormal random field.
The red line of \Cref{fig:mean_err_lognorm} shows the error decay for the standard auxiliary problem with periodic BCs and no oversampling, while the blue line displays the decay of the error for the modified elliptic method. In the first case, the only source of error is due to the BCs, while in the second case the error consists of two contributions: the boundary error which converges exponentially and is more visible in the range ($R<100$), and the systematic error dominating for larger values of $R$.

Next, we show the decay of the boundary error. The boundary error is defined as the difference $\E{ a^{0,R,L,T} - a^{0,L,T}}$ and it is controlled by an exponentially decaying deterministic upper bound. The values of $a^{0,L,T}$ are not directly accessible (they involve the solution of the cell problem over the infinite domain), so we approximate them by $a^{0,R_{max},L,T}\approx a^{0,L,T}$, with $R_{max}=500$. Additionally, we average $a^{0,R,L,T}$ and $a^{0,R_{max},L,T}$ over $N=1000$ i.i.d. samples. The exponentially decaying (in $R$) error between the empirical averages $\bar{a}^{0,R,L,T,N}$ and $\bar{a}^{0,R_{max},L,T,N}$ is depicted in Figure \ref{fig:bd_err_lognorm}.

	\subsection{Two dimensional lognormal field with exponential covariance}
\begin{figure}[h]
	\centering
	\subfloat[Logarithmic colour scale.]{
		\includegraphics[width = .45\textwidth]{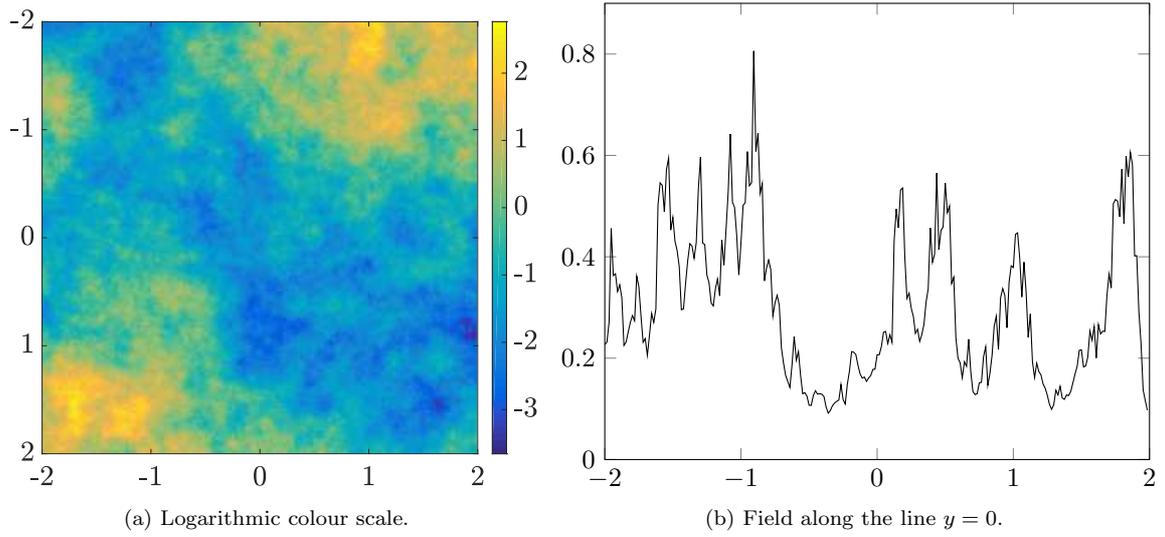}}
	\subfloat[Field along the line $y=0$.]{
%
%
\definecolor{mycolor1}{rgb}{0.00000,0.44700,0.74100}%
\begin{tikzpicture}

\begin{axis}[%
width=.45\textwidth,
height=.38\textwidth,
at={(0.521in,0.521in)},
scale only axis,
xmin=-2,
xmax=2,
ymin=0,
ymax=0.9,
axis background/.style={fill=white}
]
\addplot [color=black, forget plot]
  table[row sep=crcr]{%
-2	0.226804908217253\\
-1.984375	0.231583471983433\\
-1.96875	0.270267559948837\\
-1.953125	0.455886588133568\\
-1.9375	0.363012988064627\\
-1.921875	0.366482844458746\\
-1.90625	0.331754707172448\\
-1.890625	0.344694314716477\\
-1.875	0.319339565132577\\
-1.859375	0.225011614430954\\
-1.84375	0.231797135434793\\
-1.828125	0.251637028453652\\
-1.8125	0.270646394966996\\
-1.796875	0.284080602390156\\
-1.78125	0.273848308555158\\
-1.765625	0.361734248494757\\
-1.75	0.340973453491478\\
-1.734375	0.283407932800418\\
-1.71875	0.23123674899591\\
-1.703125	0.238728647270387\\
-1.6875	0.205755175346551\\
-1.671875	0.245517312251829\\
-1.65625	0.288285225403512\\
-1.640625	0.270470588071686\\
-1.625	0.296904796732886\\
-1.609375	0.502345180673614\\
-1.59375	0.545882241453442\\
-1.578125	0.54094149213932\\
-1.5625	0.48875454617461\\
-1.546875	0.572964888360203\\
-1.53125	0.594739217749823\\
-1.515625	0.452832636045208\\
-1.5	0.477696421267048\\
-1.484375	0.430660403387594\\
-1.46875	0.412454042045463\\
-1.453125	0.380756685245903\\
-1.4375	0.295556380467616\\
-1.421875	0.29725056721281\\
-1.40625	0.348595620335784\\
-1.390625	0.392963625084482\\
-1.375	0.425952049493404\\
-1.359375	0.421125778337285\\
-1.34375	0.39728621297603\\
-1.328125	0.427971106897858\\
-1.3125	0.534771045353181\\
-1.296875	0.596422171119955\\
-1.28125	0.427286001966516\\
-1.265625	0.42398921184038\\
-1.25	0.36402928671607\\
-1.234375	0.348897945414319\\
-1.21875	0.308881585839553\\
-1.203125	0.303160143228872\\
-1.1875	0.337235425450726\\
-1.171875	0.356069182218024\\
-1.15625	0.32305939684429\\
-1.140625	0.433349800789938\\
-1.125	0.383696597166648\\
-1.109375	0.454675041913386\\
-1.09375	0.52797397025578\\
-1.078125	0.642016679924623\\
-1.0625	0.508007789762642\\
-1.046875	0.497204226493857\\
-1.03125	0.442373517980429\\
-1.015625	0.364360847132924\\
-1	0.42661283455285\\
-0.984375	0.501374126797354\\
-0.96875	0.506797484332334\\
-0.953125	0.607479740391746\\
-0.9375	0.540390616008363\\
-0.921875	0.546311123337605\\
-0.90625	0.806089174046627\\
-0.890625	0.607941216976684\\
-0.875	0.643923614490098\\
-0.859375	0.524824360625146\\
-0.84375	0.544520616828272\\
-0.828125	0.352473619763497\\
-0.8125	0.37764455477608\\
-0.796875	0.393514571040228\\
-0.78125	0.375084682871651\\
-0.765625	0.282953913229792\\
-0.75	0.309145331205909\\
-0.734375	0.323923123465859\\
-0.71875	0.304579375868534\\
-0.703125	0.220413226095229\\
-0.6875	0.190126656721566\\
-0.671875	0.165685555690975\\
-0.65625	0.153843019042723\\
-0.640625	0.142181746231374\\
-0.625	0.187691252714021\\
-0.609375	0.242852117586591\\
-0.59375	0.196029902444039\\
-0.578125	0.212923259780961\\
-0.5625	0.166333123646307\\
-0.546875	0.129852538119597\\
-0.53125	0.13192438854616\\
-0.515625	0.125114166401867\\
-0.5	0.107142987727294\\
-0.484375	0.106880225203966\\
-0.46875	0.126310361589771\\
-0.453125	0.135851663480087\\
-0.4375	0.129048364106979\\
-0.421875	0.129850257263415\\
-0.40625	0.128910668096823\\
-0.390625	0.124277172910217\\
-0.375	0.103767035318827\\
-0.359375	0.0915113946798082\\
-0.34375	0.0966760128487246\\
-0.328125	0.106288313368026\\
-0.3125	0.111608183004132\\
-0.296875	0.114841371924313\\
-0.28125	0.117424264336047\\
-0.265625	0.147962307479903\\
-0.25	0.118329546089696\\
-0.234375	0.109524408836322\\
-0.21875	0.148878339513503\\
-0.203125	0.174703440535587\\
-0.1875	0.212772996969214\\
-0.171875	0.211763058803463\\
-0.15625	0.207099463575663\\
-0.140625	0.185878848315712\\
-0.125	0.168908302868948\\
-0.109375	0.160736412825977\\
-0.09375	0.162825072516667\\
-0.078125	0.154181215604587\\
-0.0625	0.160889016518339\\
-0.046875	0.165836128746121\\
-0.03125	0.177815256008246\\
-0.015625	0.178432883034125\\
0	0.206521393099504\\
0.015625	0.206745505078278\\
0.03125	0.22286937318514\\
0.046875	0.251069329496373\\
0.0625	0.252746163517276\\
0.078125	0.228924485386795\\
0.09375	0.246087219268799\\
0.109375	0.231193043032147\\
0.125	0.432411464566671\\
0.140625	0.4943356810652\\
0.15625	0.457369127533324\\
0.171875	0.531613762311643\\
0.1875	0.53566026492377\\
0.203125	0.400385972257707\\
0.21875	0.317930868618144\\
0.234375	0.328156881118058\\
0.25	0.299186009129554\\
0.265625	0.28266364137643\\
0.28125	0.249008349569974\\
0.296875	0.256959556959773\\
0.3125	0.233131235720642\\
0.328125	0.246793255608291\\
0.34375	0.266466885295266\\
0.359375	0.310613787004008\\
0.375	0.337116833670521\\
0.390625	0.45662847822879\\
0.40625	0.389367779153998\\
0.421875	0.405087080206776\\
0.4375	0.565260098176756\\
0.453125	0.414255255728569\\
0.46875	0.452050794436921\\
0.484375	0.458871530870658\\
0.5	0.545150568937177\\
0.515625	0.487388163144908\\
0.53125	0.50004875800583\\
0.546875	0.347591478363192\\
0.5625	0.359616171427818\\
0.578125	0.239894358208343\\
0.59375	0.199696598690483\\
0.609375	0.190573060267604\\
0.625	0.161426841804297\\
0.640625	0.192840323078351\\
0.65625	0.183959708066898\\
0.671875	0.237318823461938\\
0.6875	0.172753606280096\\
0.703125	0.134198356143303\\
0.71875	0.128070171846368\\
0.734375	0.129901604268527\\
0.75	0.151525219985701\\
0.765625	0.152178300406665\\
0.78125	0.203228872688286\\
0.796875	0.222423935155814\\
0.8125	0.146575513435626\\
0.828125	0.170306624235871\\
0.84375	0.205670856748272\\
0.859375	0.263664304927401\\
0.875	0.282772451924431\\
0.890625	0.217073191650992\\
0.90625	0.318990559302239\\
0.921875	0.337322504509878\\
0.9375	0.323879697629703\\
0.953125	0.260017366541218\\
0.96875	0.349895083908863\\
0.984375	0.382294998126295\\
1	0.378967634883279\\
1.015625	0.444233018952112\\
1.03125	0.447236517380785\\
1.046875	0.397607622830811\\
1.0625	0.320686216602831\\
1.078125	0.389531965190276\\
1.09375	0.331307291008197\\
1.109375	0.241266140180438\\
1.125	0.287128980111302\\
1.140625	0.19291236645955\\
1.15625	0.16179476374838\\
1.171875	0.190004557366495\\
1.1875	0.173481517692251\\
1.203125	0.167011837760213\\
1.21875	0.148884326856648\\
1.234375	0.140524124327617\\
1.25	0.126641336770461\\
1.265625	0.108770641033369\\
1.28125	0.0991240129946236\\
1.296875	0.108017758019971\\
1.3125	0.1368200657006\\
1.328125	0.127598012091389\\
1.34375	0.144413151895326\\
1.359375	0.121444830630922\\
1.375	0.118709404293293\\
1.390625	0.127483786113274\\
1.40625	0.126133668074811\\
1.421875	0.134461302359526\\
1.4375	0.151227866559813\\
1.453125	0.168524459895234\\
1.46875	0.183666268454961\\
1.484375	0.215454620926394\\
1.5	0.216359386480652\\
1.515625	0.182759066981011\\
1.53125	0.18455580538548\\
1.546875	0.202098523164562\\
1.5625	0.253808966578749\\
1.578125	0.235414078550249\\
1.59375	0.262395267775531\\
1.609375	0.200398397987654\\
1.625	0.266236446009645\\
1.640625	0.247871807713202\\
1.65625	0.253235433003429\\
1.671875	0.252696663577593\\
1.6875	0.324195385364052\\
1.703125	0.366025778054801\\
1.71875	0.338658453950427\\
1.734375	0.504878791059332\\
1.75	0.512972646358047\\
1.765625	0.509603179693101\\
1.78125	0.479088044456983\\
1.796875	0.572994147017084\\
1.8125	0.464958550590462\\
1.828125	0.598183954435274\\
1.84375	0.557691654726511\\
1.859375	0.605968927498204\\
1.875	0.580706673672172\\
1.890625	0.401745320301612\\
1.90625	0.402053604269982\\
1.921875	0.291648200348788\\
1.9375	0.229965642197693\\
1.953125	0.137258436199609\\
1.96875	0.116351543896875\\
1.984375	0.0974372202461841\\
};
\end{axis}
\end{tikzpicture}%
}
	\caption{Lognormal random field with exponential covariance function.}
	\label{fig: 2D lognormal field}
\end{figure}

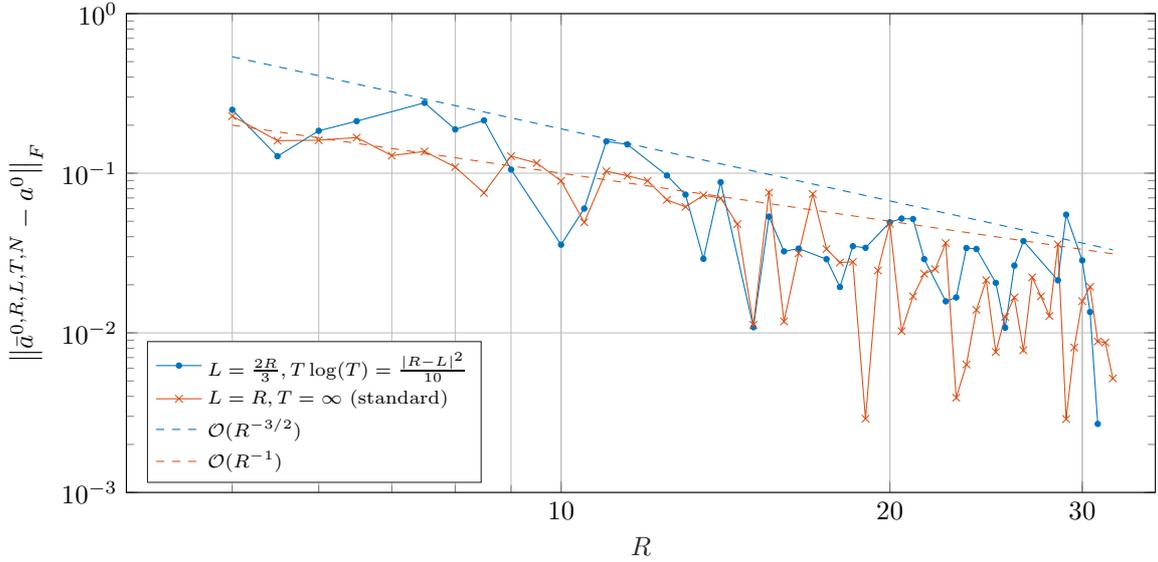
\begin{figure}[h!]
	\centering	
%
%
\definecolor{mycolor1}{rgb}{0.00000,0.44700,0.74100}%
\definecolor{mycolor2}{rgb}{0.85000,0.32500,0.09800}%
\definecolor{mycolor3}{rgb}{0.49400,0.18400,0.55600}%
\begin{tikzpicture}

\begin{axis}[%
width=.85\linewidth,
height=2.5in,
scale only axis,
xmode=log,
xmin=4,
xmax=35,
xminorticks=true,
xlabel style={font=\color{white!15!black}},
xlabel={$R$},
xtick={1, 10, 20, 30},
xticklabels={1, 10, 20, 30},
minor xtick={1,2,3,4,5,6,7,8,9,10,20,30},
ymode=log,
ymin=.001,
ymax=1,
yminorticks=true,
axis background/.style={fill=white},
title style={font=\bfseries},
ylabel={$\norm{ \bar{a}^{0,R,L,T,N} - a^{0}}_F$},
xmajorgrids,
xminorgrids,
ymajorgrids,
legend style={font=\scriptsize, legend cell align=left, align=left, draw=white!15!black,at={(0.02,0.02)},anchor=south west }
]

\addplot [color=mycolor1, mark=*, mark size=1pt, mark options={solid, mycolor1}]
  table[row sep=crcr]{%
5	0.249641423586473\\
5.5	0.127813296793226\\
6	0.184362696595248\\
6.5	0.211945485997345\\
7.5	0.276179508014084\\
8	0.188150843169365\\
8.5	0.214360625891067\\
9	0.10542355373622\\
10	0.03566713735546\\
10.5	0.0600511867127295\\
11	0.15820099351577\\
11.5	0.151607050463648\\
12.5	0.0967827967377502\\
13	0.0734339699738979\\
13.5	0.0290965494414664\\
14	0.0878533247654354\\
15	0.010839584397971\\
15.5	0.0534647775386522\\
16	0.0323907372023164\\
16.5	0.0337033641275292\\
17.5	0.0289394347771064\\
18	0.0193691727383057\\
18.5	0.0348981982263486\\
19	0.034105858045153\\
20	0.0492793394870546\\
20.5	0.0519865975495229\\
21	0.0516616170679634\\
21.5	0.0289682398230567\\
22.5	0.015757530833393\\
23	0.0166747379437859\\
23.5	0.034013870112555\\
24	0.0335472926028613\\
25	0.0205824802552731\\
25.5	0.0107755054165744\\
26	0.0263954387098743\\
26.5	0.037579620700242\\
28.5	0.0213249528337653\\
29	0.0550138980026853\\
30	0.0284380330475002\\
30.5	0.0135069776412551\\
31	0.00268731698994223\\
};
\addlegendentry{$L=\frac{2R}{3}, T\log(T) = \frac{\abs{R-L}^2}{10}$}

\addplot [color=mycolor2, mark=x, mark options={solid, mycolor2}]
  table[row sep=crcr]{%
5	0.2279656483464\\
5.5	0.159975202711929\\
6	0.160983355689051\\
6.5	0.167299708505859\\
7	0.129014718942545\\
7.5	0.136931034137876\\
8	0.109373958265985\\
8.5	0.0752389495974278\\
9	0.127924470500264\\
9.5	0.11585675525826\\
10	0.089566233979057\\
10.5	0.0492555678750969\\
11	0.103062304329046\\
11.5	0.0963497627022967\\
12	0.0896316118707314\\
12.5	0.0681057158786466\\
13	0.0615346647307144\\
13.5	0.0727674906346586\\
14	0.0697284141939716\\
14.5	0.0479660085922424\\
15	0.0112142941110742\\
15.5	0.0756754093194552\\
16	0.011808028120563\\
16.5	0.0315122640286781\\
17	0.0741304884004592\\
17.5	0.0335095336305544\\
18	0.0275806673069093\\
18.5	0.027842783681943\\
19	0.00289691041137363\\
19.5	0.0245578022960984\\
20	0.0482502998190108\\
20.5	0.0102412527875059\\
21	0.0169295171372911\\
21.5	0.0234518487889467\\
22	0.0250621821194206\\
22.5	0.0365013244940011\\
23	0.00392342536462712\\
23.5	0.00633650795060869\\
24	0.0139307357599865\\
24.5	0.0214010827287769\\
25	0.00758861526529982\\
25.5	0.0126001242492804\\
26	0.0166786400621708\\
26.5	0.00777569304637885\\
27	0.022272070704311\\
27.5	0.0169470493962706\\
28	0.0127482830688092\\
28.5	0.035852043599155\\
29	0.0028843619698292\\
29.5	0.00809537532250028\\
30	0.0157917151975397\\
30.5	0.0194303325907455\\
31	0.00883977482491628\\
31.5	0.00870043901081582\\
32	0.00518236646113233\\
};
\addlegendentry{$L=R, T = \infty$ (standard)}

\addplot [color=mycolor1, dashed]
  table[row sep=crcr]{%
5	0.53665631459995\\
5.5	0.465165199321332\\
6	0.408248290463863\\
6.5	0.362060557178186\\
7	0.323969548293623\\
7.5	0.292118697336089\\
8	0.265165042944955\\
8.5	0.242115649612954\\
9	0.222222222222222\\
9.5	0.204911268796753\\
10	0.189736659610103\\
10.5	0.176346685670962\\
11	0.164460733406053\\
11.5	0.153852475987719\\
12	0.144337567297406\\
12.5	0.135764501987817\\
13	0.128007737590437\\
13.5	0.120962456433737\\
14	0.114540532248182\\
14.5	0.108667392711495\\
15	0.103279555898864\\
15.5	0.0983226789678894\\
16	0.09375\\
16.5	0.0895210843486056\\
17	0.0856008088363528\\
17.5	0.0819585332115013\\
18	0.0785674201318386\\
18.5	0.0754038737797747\\
19	0.0724470738538616\\
19.5	0.0696785867165905\\
20	0.0670820393249937\\
20.5	0.0646428445316175\\
21	0.062347968638855\\
21.5	0.0601857338809889\\
22	0.0581456499151665\\
22.5	0.0562182695141045\\
23	0.0543950645366282\\
23.5	0.0526683189578958\\
24	0.0510310363079829\\
24.5	0.0494768593250062\\
25	0.048\\
25.5	0.04659517849302\\
26	0.0452575696472732\\
26.5	0.0439827560312723\\
27	0.0427666866066389\\
27.5	0.0416056402579895\\
28	0.0404961935367029\\
28.5	0.0394351920665976\\
29	0.0384197251400797\\
29.5	0.0374471031009487\\
30	0.0365148371670111\\
30.5	0.0356206213937778\\
31	0.0347623165213113\\
31.5	0.0339379354809428\\
32	0.0331456303681194\\
};
\addlegendentry{$\mathcal{O}(R^{-3/2})$}

\addplot [color=mycolor2, dashed]
table[row sep=crcr]{%
	5	0.2\\
	5.5	0.181818181818182\\
	6	0.166666666666667\\
	6.5	0.153846153846154\\
	7	0.142857142857143\\
	7.5	0.133333333333333\\
	8	0.125\\
	8.5	0.117647058823529\\
	9	0.111111111111111\\
	9.5	0.105263157894737\\
	10	0.1\\
	10.5	0.0952380952380952\\
	11	0.0909090909090909\\
	11.5	0.0869565217391304\\
	12	0.0833333333333333\\
	12.5	0.08\\
	13	0.0769230769230769\\
	13.5	0.0740740740740741\\
	14	0.0714285714285714\\
	14.5	0.0689655172413793\\
	15	0.0666666666666667\\
	15.5	0.0645161290322581\\
	16	0.0625\\
	16.5	0.0606060606060606\\
	17	0.0588235294117647\\
	17.5	0.0571428571428571\\
	18	0.0555555555555556\\
	18.5	0.0540540540540541\\
	19	0.0526315789473684\\
	19.5	0.0512820512820513\\
	20	0.05\\
	20.5	0.0487804878048781\\
	21	0.0476190476190476\\
	21.5	0.0465116279069767\\
	22	0.0454545454545455\\
	22.5	0.0444444444444444\\
	23	0.0434782608695652\\
	23.5	0.0425531914893617\\
	24	0.0416666666666667\\
	24.5	0.0408163265306122\\
	25	0.04\\
	25.5	0.0392156862745098\\
	26	0.0384615384615385\\
	26.5	0.0377358490566038\\
	27	0.037037037037037\\
	27.5	0.0363636363636364\\
	28	0.0357142857142857\\
	28.5	0.0350877192982456\\
	29	0.0344827586206897\\
	29.5	0.0338983050847458\\
	30	0.0333333333333333\\
	30.5	0.0327868852459016\\
	31	0.032258064516129\\
	31.5	0.0317460317460317\\
	32	0.03125\\
};
\addlegendentry{$\mathcal{O}(R^{-1})$}

\end{axis}
\end{tikzpicture}%
	\caption{Comparison of the error in mean of the standard and modified elliptic methods for a two-dimensional lognormal random field with exponential covariance function.}
	\label{fig: lognormal 2d err mean}
\end{figure}

	We will now study the convergence for a two dimensional lognormal random field with exponential covariance function, such as the one depicted in \Cref{fig: 2D lognormal field}. 
The field is sampled by generating a Gaussian random field over the uniform grid
$$
\set{ \left(x_i,y_j\right)\in K_R: x_i=ih-\frac{R}{2}, y_j=jh-\frac{R}{2}},
$$
coinciding with the set of vertices of the structured mesh of stepsize $h=2^{-5}$ on $K_R$.
The correctors are computed by the finite element method with $\mathbb{P}_1$-elements, and the right-hand side is calculated by the Krylov subspace method with up to 2000 basis elements. 
The average is approximated by drawing $N=200$ i.i.d. samples of the lognormal field.
The convergence behaviour of the mean error is pictured in \Cref{fig: lognormal 2d err mean} for both the method discussed in this work and the truncated domain approach discussed in \Cref{subsec: naive elliptic}.
%
In this case we notice that the convergence rates improve to $1$ and $3/2$ for, respectively, the standard and the modified elliptic approaches. 

\section{Conclusion}
In this work we addressed the problem of estimating the systematic error for the modified elliptic model, defined as the difference between the true effective coefficient $a^0$ and its approximation by the corrector problem \eqref{eq: modified_Elliptic} over the infinite domain $\R^d$. 
By exploiting the time decay properties of solutions of linear parabolic equations, we found that the systematic error scales as $T^{-d/2}$, where $T$ is the final time. The parabolic solution, evaluated at time $T$, enters as a source term into \eqref{eq: modified_Elliptic}.
We have also prove the existence of a corrector $\nabla\chi_{T}\in \mathcal{L}^2_{pot}$ using the very same time decay properties.
The numerical experiments demonstrate that the modified elliptic approach outperforms the standard one in terms of decay rate of the error as the cell size $R$ grows. This is achieved upon choosing $T$ such that 
$$
c\abs{R-L} \le T \le C \abs{R-L}^2,
$$
for some $c,C>0$. Moreover, the theoretical bound of \Cref{thm: systematic error} is verified. The two dimensional numerical results presented here also verify the improved convergence rates in higher dimensions. 


\section{Acknowledgments}
The authors are grateful to J.~C.~Mourrat for precious help and suggestions. This work is partially supported by the Swiss National Science Foundation grant No. 200020\_172710.

\bibliographystyle{plain}
\bibliography{./references}

\end{document}